\documentclass[12pt]{article}
\usepackage{amsmath,amsthm,amssymb}
\usepackage{epic}
\usepackage{eepic}
\usepackage{float}
\usepackage{graphicx}
\makeatletter
\renewcommand{\@makecaption}[2]{%
  \vskip\abovecaptionskip
  \sbox\@tempboxa{#1 #2}%
  \ifdim \wd\@tempboxa >\hsize
    #1: #2\par
  \else
    \global \@minipagefalse
    \hb@xt@\hsize{\hfil\box\@tempboxa\hfil}%
  \fi
  \vskip\belowcaptionskip}
\makeatother
\theoremstyle{definition}
\newtheorem{definition}{Definition}[section]
\newtheorem{example}[definition]{Example}
\theoremstyle{plain}
\newtheorem{theorem}[definition]{Theorem}
\newtheorem{corollary}[definition]{Corollary}
\newtheorem{fact}[definition]{Fact}
\theoremstyle{definition}
\newtheorem{defin}{Definition}[subsection]
\theoremstyle{remark}
\newtheorem{rem}[defin]{Remark}
\newtheorem{Rem}[definition]{Remark}
\theoremstyle{plain}
\newtheorem{lemm}[defin]{Lemma}
\theoremstyle{plain}
\newtheorem{claimsubsec}{Claim}

\newtheorem{prop}[definition]{Proposition}
\newtheorem{lemma}{Lemma}

\theoremstyle{remark}
\newtheorem*{remark}{Remark}
\def   \set#1#2{\{{#1}:{#2}\}}
\def \itm#1{\newline\noindent{\rm{#1}}\enspace}

\newcommand{\rrank}{\mathbb{R}\text{-}\operatorname{rank}}

\newcommand{\Aut}{\operatorname{{Aut}}}
\newcommand{\diag}{\operatorname{diag}}
\newcommand{\rarrowsim}{\smash{\mathop{\,\longrightarrow\,}\limits
  ^{\lower1.5pt\hbox{$\scriptstyle\sim$}}}}
\providecommand{\bysame}{\makebox[3em]{\hrulefill}\thinspace}
\newenvironment{eq-text}{\par
\refstepcounter{equation}
\noindent
{\upshape(\theequation)}\quad}
{\par
\noindent \ignorespacesafterend}
\title{Visible actions on symmetric spaces}
\author{Toshiyuki KOBAYASHI\\[\smallskipamount]
Research Institute for Mathematical Sciences\\
Kyoto University}
\date{}
\begin{document}
\maketitle

\numberwithin{equation}{subsection}
\numberwithin{table}{subsection}

\begin{abstract}
A visible action on a complex manifold is a holomorphic action
that admits a $J$-transversal totally real submanifold $S$.
It is said to be strongly visible if there exists an orbit-preserving
anti-holomorphic diffeomorphism 
$\sigma$ such that $\sigma |_S = \mathrm{id}$.

In this paper,
we prove that for any Hermitian symmetric space $D = G/K$
the action of any symmetric subgroup $H$
is strongly visible.
The proof is carried out by finding explicitly an orbit-preserving
anti-holomorphic involution
and a totally real submanifold $S$.

Our geometric results provide a uniform proof of various 
multiplicity-free theorems of irreducible highest weight modules 
when restricted to reductive symmetric pairs,
for both classical and exceptional cases,
for both finite and infinite dimensional cases,
and for both discrete and continuous spectra.
\end{abstract}



\setcounter{tocdepth}{1}
\tableofcontents

\section{Introduction and main results}
Suppose a Lie group $H$ acts holomorphically on a connected complex
manifold $D$ with complex structure $J$.
We recall from
\cite[Definition~2.3]{Acta2004}:
\begin{definition}
The action is \textbf{visible} if there exist a (non-empty)
$H$-invariant open subset $D'$ of $D$ and a totally real submanifold
$S$ of $D'$ satisfying the following two conditions:
\addtocounter{subsection}{1}
\begin{align}
& \text{$S$ meets every $H$-orbit in $D'$,}
\label{eqn:1.1.1}
\\
& \text{$J_x (T_x S) \subset T_x (H \cdot x)$ for all $x \in S$
        \quad($J$-transversality).}
\label{eqn:1.1.2}
\end{align}
\end{definition}

Obviously, a transitive action is visible.
Conversely, a visible action requires the existence of an $H$-orbit
whose dimension is at least half the real dimension of $D$.  Further,
in \cite[Definition~3.3.1]{RIMS}, we have introduced:
\begin{definition}\label{def:S}
The action 
is \textbf{strongly visible}
if there exist an $H$-invariant open subset $D'$ of $D$,
a submanifold $S$ (\textbf{slice}) of $D'$,
and an anti-holomorphic diffeomorphism $\sigma$ of $D'$
satisfying the following three conditions:
\addtocounter{subsection}{1}
\setcounter{equation}{0}
\begin{align}
& \text{$S$ meets every $H$-orbit in $D'$,}
\label{eqn:1.2.1}
\\
& \text{$\sigma|_S = \mathrm{id}$,}
\label{eqn:1.2.2}
\\
& \text{$\sigma$ preserves each $H$-orbit in $D'$.}
\label{eqn:1.2.3}
\end{align}
\end{definition}

A strongly visible action is visible 
(see \cite[Theorem~4]{RIMS}).
The concept of strongly visible actions
is used as
 a crucial assumption on base spaces $D$ for the
propagation theorem of 
multiplicity-free property from fibers to
spaces of sections of equivariant holomorphic vector
bundles over $D$ 
(see \cite{mfbdle}).

The aim of this article is to give a systematic study of strongly
visible actions on symmetric spaces. 
In a previous paper \cite[Theorem 11]{RIMS},
we have discussed the case where $D$ is a complex symmetric space
$G_{\mathbb{C}}/K_{\mathbb{C}}$:
\begin{fact}
\label{fact:cpxsymm}
Suppose $G/K$ is a Riemannian symmetric space,
and $G_{\mathbb{C}}/K_{\mathbb{C}}$ is its complexification.
Then, the $G$-action on the complex symmetric space
$G_{\mathbb{C}}/K_{\mathbb{C}}$ is strongly visible.
\end{fact}

In this article, our focus is on the case where $D$ is 
a Hermitian symmetric space.
A typical example is:
\begin{example}
\label{ex:SLtwo}
Let $G = KAN$ be the Iwasawa decomposition of
$G := SL(2,\mathbb{R})$.
Then, all of the actions of $K, A$ and $N$ on the Hermitian symmetric
space 
$D := G/K$
are strongly visible,
as one can see easily from the following figures where $D$ is realized
as the Poincar\'{e} disk and the dotted lines give slices:
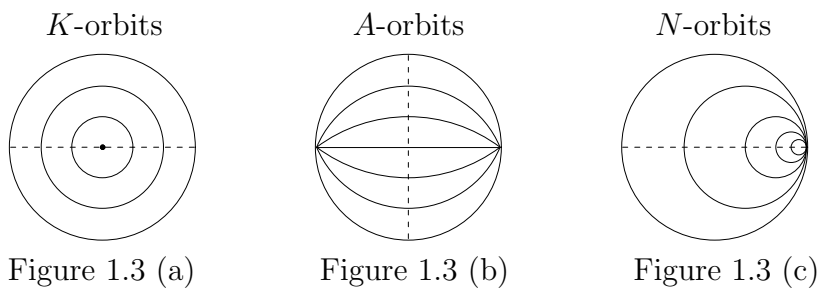
\begin{center}
$K$-orbits \hspace{5.8em}
$A$-orbits \hspace{5.8em}
$N$-orbits
\medskip

\setlength{\unitlength}{0.00053333in}
\begingroup\makeatletter\ifx\SetFigFont\undefined%
\gdef\SetFigFont#1#2#3#4#5{%
  \reset@font\fontsize{#1}{#2pt}%
  \fontfamily{#3}\fontseries{#4}\fontshape{#5}%
  \selectfont}%
\fi\endgroup%
{\renewcommand{\dashlinestretch}{30}
\begin{picture}(7840,1853)(0,-10)
\dashline{60.000}(6020,919)(7820,919)
\dashline{60.000}(3920,1819)(3920,19)
\dashline{60.000}(20,919)(1820,919)
\path(3020,919)(4820,919)
\put(923,918){\blacken\ellipse{48}{48}}
\put(923,918){\ellipse{48}{48}}
\put(7520,919){\ellipse{600}{600}}
\put(7670,919){\ellipse{300}{300}}
\put(7220,919){\ellipse{1200}{1200}}
\put(7745,919){\ellipse{150}{150}}
\put(6920,919){\ellipse{1824}{1824}}
\put(3920,919){\ellipse{1824}{1824}}
\put(920,919){\ellipse{1824}{1824}}
\put(920,919){\ellipse{1200}{1200}}
\put(920,919){\ellipse{600}{600}}
\put(3920.000,1294.000){\arc{1950.000}{0.3948}{2.7468}}
\put(3920.000,2119.000){\arc{3000.000}{0.9273}{2.2143}}
\put(3920.000,544.000){\arc{1950.000}{3.5364}{5.8884}}
\put(3920.000,-281.000){\arc{3000.000}{4.0689}{5.3559}}
\end{picture}
}
%

\centerline{Figure 1.3 (a)\hspace{4em}
Figure 1.3 (b)\hspace{4em}
Figure 1.3 (c)}
\end{center}
\end{example}

In this case, both $(G,K)$ and $(G,A)$ form symmetric pairs,
while $N$ is a maximal unipotent subgroup of $G$.

These three examples are generalized to the following
Theorems~\ref{thm:Hvisible} and \ref{thm:Nvisible},
which are our main results of this paper.
\begin{theorem}
\label{thm:Hvisible}
Suppose $G$ is a semisimple Lie group such that
 $D := G/K$ is a Hermitian symmetric space.
Then, for any symmetric pair $(G,H)$,
 the $H$-action on $D$ is strongly visible.
\end{theorem}

\begin{example}[the Siegel upper half space]
\label{ex:Hvisible}
Let $G/K = Sp(n,\mathbb{R})/U(n)$. 
Then, the action of the subgroup $H$  is strongly visible if
$H = GL(n,\mathbb{R})$,
$U(p,q), Sp(p,\mathbb{R}) \times Sp(q,\mathbb{R})$
$(p+q = n)$,
or
$Sp(\frac{n}{2}, \mathbb{C})$ $(n:$ even$)$.
\end{example}

The pair $(G \times G, \diag(G))$ is a classic example of symmetric pairs, 
where
$\diag(G) := \{(g,g): g \in G\}$.
Thus, Theorem~\ref{thm:Hvisible} also includes:
\begin{theorem}
\label{thm:dvisible}
Let $D_1, D_2$ be two Hermitian symmetric spaces of a compact simple
 Lie group
$G_U$.
Then the diagonal action of $G_U$ on $D_1 \times D_2$ is strongly
visible.
\end{theorem}

\begin{example}
\label{ex:dvisible}
The diagonal action of $SU(n)$ on the direct product of two Grassmann
varieties 
$Gr_p (\mathbb{C}^n) \times Gr_k (\mathbb{C}^n)$
$(1 \le p, k \le n)$
is strongly visible.
\end{example}

\begin{theorem}
\label{thm:Gvisible}
Let $D$ be a non-compact irreducible Hermitian symmetric space $G/K$,
and $\overline{D}$ denote the Hermitian symmetric space equipped with
the opposite complex structure.

\begin{enumerate}
    \renewcommand{\labelenumi}{\upshape \theenumi)}
\item 
The diagonal action of $G$ on $D \times D$ is strongly visible.
\item 
The diagonal action of $G$ on $D \times \overline{D}$ is strongly visible. 
\end{enumerate}
\end{theorem}

The third case of Example~\ref{ex:SLtwo} is generalized as follows: 

\begin{theorem}
\label{thm:Nvisible}
Suppose $D = G/K$ is a Hermitian symmetric space without compact
factor.
If $N$ is a maximal unipotent subgroup of $G$,
then the $N$-action on $D$ is strongly visible.
\end{theorem}

This paper is organized as follows:
In Section~\ref{sec:2},
we translate geometric conditions of (strongly) visible actions
into an algebraic language by using the structural theory of semisimple
symmetric pairs.
The proof of our main result,
Theorem \ref{thm:Hvisible}, is given in Section~\ref{sec:3}
(non-compact case) and in Section~\ref{sec:4} (compact case),
and that of Theorem \ref{thm:Nvisible} is given in
Section~\ref{sec:5}. 

As Fact \ref{fact:cpxsymm} gives a new proof (see \cite{mfbdle})
of the Cartan--Gelfand
theorem that the Plancherel formula for a Riemannian symmetric space
is multiplicity-free (induction of representations),
our geometric results here give a number of multiplicity-free theorems,
in particular, in branching problems (restriction of representations)
for both finite and infinite dimensional representations
and for discrete and continuous spectra.
Such applications to representation theory are discussed
in Section~\ref{sec:6}.

Concepts related to visible actions on complex manifolds are polar
actions on Riemannian manifolds and coisotropic actions on symplectic
manifolds. 
Since Hermitian symmetric spaces $D$ are K\"{a}hler,
we can compare these three concepts for $D$.
Some comments on this are given in Section~\ref{sec:cpv}.

\section{Preliminary results}
\label{sec:2}
This section provides sufficient conditions by means of Lie algebras
for the geometric conditions \eqref{eqn:1.1.1} -- \eqref{eqn:1.2.3}
for strongly visible actions in the setting where $D$ is a Hermitian
symmetric space.

\subsection{Semisimple symmetric pairs}
\label{subsec:3.1}
Let $G$ be a semisimple Lie group with Lie algebra $\mathfrak{g}$.
Suppose that $\tau$ is an involutive automorphism of $G$.
We write  $G_0^\tau$ for
the identity component of
$G^\tau := \{g \in G : \tau g = g \}$.
The pair $(G,H)$ 
(or the pair $(\mathfrak{g}, \mathfrak{h})$ of their Lie algebras)
is called a \textit{(semisimple) symmetric pair}
if a subgroup $H$ satisfies $G_0^\tau \subset H \subset G^\tau$.
Unless otherwise mentioned,
we shall take $H$ to be $G_0^\tau$
because Theorem~\ref{thm:Hvisible} follows readily from this case.
We  use the same letter $\tau$ to denote its differential,
and set
$$
\mathfrak{g}^{\pm \tau} := \set{Y \in \mathfrak{g}}{\tau Y = \pm Y}.
$$
Then,
$\mathfrak{g}^\tau$ is the Lie algebra of $H$.
Since $\tau^2 = \mathrm{id}$,
 we have a direct sum decomposition $\mathfrak{g} = \mathfrak{g}^\tau + \mathfrak{g}^{-\tau}$.

It is known that there exists
 a Cartan involution $\theta$ of $G$ commuting with $\tau$.
Take such $\theta$, and we write 
 $K := G^\theta =\set{g \in G}{\theta g = g}$.
Then we have a Cartan decomposition
  $\mathfrak{g} = \mathfrak{k}+ \mathfrak{p} \equiv \mathfrak{g}^\theta + \mathfrak{g}^{-\theta}$.
We shall allow $G$ to be non-linear,
and therefore
 $K$ is not necessarily compact.
The \textit{real rank} of $\mathfrak{g}$,
 denoted by $\rrank \mathfrak{g}$,
 is defined to be the dimension of a maximal abelian subspace of
 $\mathfrak{g}^{-\theta}$.

As $(\tau \theta)^2 = \mathrm{id}$,
the pair $(\mathfrak{g},\mathfrak{g}^{\tau\theta})$
also forms a symmetric pair.
The Lie group 
$G^{\tau\theta} = \{ g \in G: (\tau \theta)g = g \}$
is a reductive Lie group with Cartan involution
$\theta |_{G^{\tau\theta}}$,
and its
 Lie algebra $\mathfrak{g}^{\tau\theta}$ is  reductive 
with Cartan decomposition
$$
\mathfrak{g}^{\tau\theta} = \mathfrak{g}^{\theta,\tau\theta} +
\mathfrak{g}^{-\theta,\tau\theta} = \mathfrak{g}^{\theta,\tau} +
\mathfrak{g}^{-\theta,-\tau} ,
$$
where we have used the notation 
$\mathfrak{g}^{-\theta,-\tau}$ and alike,
defined as follows:
$$
   \mathfrak{g}^{-\theta, -\tau}
   :=
   \set{Y \in \mathfrak{g}}{(-\theta) Y = (-\tau) Y = Y}.
$$
The real rank of $\mathfrak{g}^{\theta\tau}$ is referred to as the 
\textit{split rank} of the symmetric space $G/H$,
denoted by $\mathbb{R}$-rank $G/H$ or
$\mathbb{R}$-rank $\mathfrak{g}/\mathfrak{g}^\tau$. 
That is,
\begin{equation}
    \rrank \mathfrak{g}^{\tau\theta} = \rrank \mathfrak{g}/\mathfrak{g}^\tau
 .
\label{eqn:2.1.1}
\end{equation}
In particular,
 $\rrank \mathfrak{g} = \rrank \mathfrak{g}/\mathfrak{k}$ 
if we take $\tau$ to be the Cartan involution $\theta$.

\subsection{Stability of $H$-orbits}
\label{subsec:3.2}
Retain the setting in Subsection~\ref{subsec:3.1}.
Suppose furthermore that there is another
  involutive automorphism $\sigma$ of $G$ such that
$\sigma\theta = \theta\sigma$ and
$\sigma\tau = \tau\sigma$.
The commutativity of $\sigma$ and $\theta$ implies that
 the automorphism $\sigma$ 
 stabilizes $K$, and therefore induces a diffeomorphism of $G/K$,
 for which we shall use the same letter $\sigma$.
Then, $\sigma$ sends $H$-orbits on $G/K$ to $H$-orbits because
$\sigma\tau = \tau\sigma$.
However, 
$\sigma$ may permute $H$-orbits,
and may not preserve each $H$-orbit.
In this subsection,
we give a sufficient condition
 in terms of the real
 rank condition that $\sigma$
 preserves each $H$-orbit.
The conclusion of Lemma~\ref{lem:3.2} 
meets the requirements \eqref{eqn:1.2.1}--\eqref{eqn:1.2.3}
of strongly visible conditions, and thus
Lemma~\ref{lem:3.2} will play
 a key role in the proof of our main theorem.

\begin{lemma}
\label{lem:3.2}
Let $G$ be a semisimple Lie group,
$\theta$ a Cartan involution,
and $D := G/K$ the corresponding Riemannian symmetric space.
Let $\sigma$ and $\tau$ be involutive automorphisms of 
 $G$.
We set $H := G_0^\tau$.
We assume 
the following two conditions:
\begin{align}
& \text{$\sigma$, $\tau$ and $\theta$ commute with one another.}
\label{eqn:3.2.1}
\\
& \text{$\rrank \mathfrak{g}^{\tau\theta}
 = \rrank \mathfrak{g}^{\sigma, \tau\theta}$.}
\label{eqn:3.2.2}
\end{align}
Then the followings hold:

\noindent
\textup{1)}\enspace
For any $x \in D$, there exists $g \in H$ such that
$
      \sigma (x) = g \cdot x.
$
In particular,
$\sigma$ preserves each $H$-orbit on $D$.
\newline
\textup{2)}\enspace
Let $\mathfrak{a}$ be a maximal abelian subspace 
in 
\begin{equation}
\label{eqn:3.2.4}
\mathfrak{g}^{-\theta,\sigma,\tau\theta} =
   \mathfrak{g}^{-\theta, \sigma, -\tau}
   :=
   \set{Y \in \mathfrak{g}}{(-\theta) Y = \sigma Y = (-\tau) Y = Y} 
   .
\end{equation}
Then, the submanifold
$S := (\exp \mathfrak{a})K$
 meets every $H$-orbit in $D$ and 
$\sigma|_S = \mathrm{id}$. 
%
%
\end{lemma}

\begin{proof}
1)\enspace
First, let us show that if $h \in H$ then
$g := \sigma(h)h^{-1} \in H$.
In fact,
 by using $\sigma \tau = \tau \sigma$ and $\tau (h) = h$,
 we have
$$
 \tau (g) = \tau \sigma(h)\; \tau(h^{-1}) 
          = \sigma \tau (h)\; \tau(h)^{-1}
          = \sigma (h)\; h^{-1} = g  .
$$
Hence, 
 $g \in G^\tau$.
Moreover,
 since the image of the continuous map
$$
     H \to G,\ \ \ h \mapsto \sigma(h)\; h^{-1}
$$
 is connected,
 we have proved $g \in G^\tau_0 = H$.

Next, let
 $\mathfrak{a}$ be as in 2).
In light of the Cartan decomposition
$\mathfrak{g}^{\sigma,\tau\theta}
 = \mathfrak{g}^{\theta,\sigma,\tau\theta}
 + \mathfrak{g}^{-\theta,\sigma,\tau\theta} $, we have
$\dim \mathfrak{a} = \rrank \mathfrak{g}^{\sigma,\tau\theta}$.
Then,
 the assumption \eqref{eqn:3.2.2}  
 shows $\dim \mathfrak{a} = \rrank \mathfrak{g}^{\tau\theta}$.
As 
$\mathfrak{g}^{\tau\theta} = \mathfrak{g}^{\theta,\tau}
                            + \mathfrak{g}^{-\theta,-\tau}$
is a Cartan decomposition of 
$\mathfrak{g}^{\tau\theta}$, this 
means that $\mathfrak{a}$ is also a maximal abelian subspace in
$
   \mathfrak{g}^{-\theta, -\tau}
$.

We write $A$ for the analytic subgroup of $G$ with Lie algebra $\mathfrak{a}$.
Then we have a generalized Cartan decomposition (see \cite[\S 2]{xfjd})
\begin{equation}
   G = H A K .
\label{eqn:3.2.3}
\end{equation}
Fix $x \in G/K$.
Then,
 according to the decomposition \eqref{eqn:3.2.3}, 
 we find $h \in H$ and $a \in A$ such that
$$
        x = h a \cdot o ,
$$
where $o := e K \in G/K$.

On the other hand,
 we have $\sigma(a) = a$ because $\mathfrak{a} \subset
   \mathfrak{g}^{-\theta, \sigma, -\tau} \subset
 \mathfrak{g}^\sigma$.
Therefore, if we set
$g := \sigma(h)h^{-1} \in G_0^\tau$,
then
$$
    \sigma (x) = \sigma(h)\; \sigma(a) \cdot o
               = \sigma(h)\; h^{-1} h a \cdot o
               = g \cdot x  .
$$
In particular,
we have 
$\sigma(H \cdot x) = H \cdot \sigma(x)
  = H g \cdot x = H \cdot x$.
Thus, $\sigma$ preserves each
$H$-orbit on $G/K$.

2)\enspace
The submanifold $S$ meets every $H$-orbit by \eqref{eqn:3.2.3},
and $\sigma|_S = \mathrm{id}$ because
 $\mathfrak{a} \subset \mathfrak{g}^\sigma$.

\end{proof}

\subsection{Involutions on Hermitian Symmetric Space $G/K$}
\label{subsec:3.3}
This subsection gives a brief review on basic results on 
submanifolds of Hermitian
symmetric spaces.

Let $G$ be a non-compact simple Lie group $G$ with Cartan decomposition
$\mathfrak{g} = \mathfrak{k} + \mathfrak{p}$. 
$G$ is said to be of 
\textit{Hermitian type} if the center $\mathfrak{c}(\mathfrak{k})$ of
$\mathfrak{k}$ is non-trivial.
Then, it is known  that $\dim \mathfrak{c}(\mathfrak{k}) = 1$
and that there exists a characteristic element
$Z \in \mathfrak{c}(\mathfrak{k})$ such that
\begin{equation}
\label{eqn:GKPP}
\mathfrak{g}_{\mathbb{C}} := \mathfrak{g} \otimes \mathbb{C}
 = \mathfrak{k}_{\mathbb{C}} \oplus \mathfrak{p}_+ \oplus
 \mathfrak{p}_-
\end{equation}
is the eigenspace decomposition of
$\operatorname{ad}(Z)$ with eigenvalues $0, \sqrt{-1}$ and $-\sqrt{-1}$,
respectively.

Let $G_\mathbb{C}$ be a connected complex Lie group with Lie algebra
 $\mathfrak{g}_\mathbb{C}$,
 and $Q^-$ the parabolic subgroup of $G_\mathbb{C}$ 
 with Lie algebra $\mathfrak{k}_\mathbb{C} + \mathfrak{p}_-$.
Then the natural homomorphism
$G \to G_{\mathbb{C}}$ induces an open
 embedding $G/K \hookrightarrow G_\mathbb{C}/Q^-$,
 from which the complex structure on $G/K$ is induced.
This complex structure on $G/K$ is given by the left
$G$-translation of
\begin{equation}
\label{eqn:cpxstr}
J_o := \operatorname{Ad} (\exp (\frac{\pi}{2} Z)):
   T_o (G/K) \to T_o (G/K)
\end{equation}
at the origin $o = eK \in G/K$.

Suppose $\tau$ is an involutive automorphism of $G$.
We may and do assume that $\tau$ commutes with $\theta$
(by taking a conjugation by an inner automorphism if necessary).
Then $\tau$ stabilizes the Cartan decomposition
$\mathfrak{g} = \mathfrak{k} + \mathfrak{p}$,
and particularly the one dimensional subspace
$\mathfrak{c}(\mathfrak{k})$.
Since $\tau^2 = \mathrm{id}$,
we have either
\begin{align}
    \tau Z &= Z  ,
\label{eqn:2.3.1}
\\
\intertext{or}
   \tau Z &= -Z  .
\label{eqn:2.3.2}
\end{align}

It follows from the definition \eqref{eqn:cpxstr} of the complex
structure on $G/K$
that the condition~\eqref{eqn:2.3.1} has the following geometric
meaning: 
\begin{itemize}
\item[]
$\tau$ acts \textbf{holomorphically} on the Hermitian symmetric space
$G/K$,
\item[]
$G^\tau/K^\tau \hookrightarrow G/K$
defines a complex submanifold.
\end{itemize}
On the other hand, 
the condition~\eqref{eqn:2.3.2} implies
\begin{itemize}
\item[]
$\tau$ acts \textbf{anti-holomorphically} on the Hermitian symmetric
space $G/K$,
\item[]
$G^\tau/K^\tau \hookrightarrow G/K$ defines a totally real submanifold.
\end{itemize}
We say the involution $\tau$ (or the symmetric pair
$(\mathfrak{g}, \mathfrak{g}^\tau)$)
is of \textit{holomorphic type}
(respectively, \textit{anti-holomorphic type})
if $\tau$ satisfies \eqref{eqn:2.3.1}
(respectively, \eqref{eqn:2.3.2}).

The following Tables~\ref{tbl:3.3.1} and \ref{tbl:3.3.2} give
 the classification of semisimple symmetric pairs
 $(\mathfrak{g}, \mathfrak{g}^\tau)$ 
for simple Lie algebras $\mathfrak{g}$ 
 such that the pair $(\mathfrak{g}, \mathfrak{g}^\tau)$
is of holomorphic type and of anti-holomorphic type,
respectively.
Table~\ref{tbl:3.3.2} 
is equivalent to the classification of
 totally real symmetric spaces $G^\tau/K^\tau$ of a Hermitian
 symmetric space $G/K$
 (see \cite{xfo}, \cite{xjafbams}, \cite{xjafjdg}, \cite{xkobanaga}).
For later purposes, we label these symmetric spaces in the left column
of Tables~\ref{tbl:3.3.1} and \ref{tbl:3.3.2}.

\begin{table}[H]
$$
\vbox{
\offinterlineskip
\def\tablerule{\noalign{\hrule}}
\def\ct{&\cr\tablerule}
\halign{\strut#&\vrule#&
            \;\;\hfil#\hfil\hfil\;\;&\vrule#&
            \;\;\hfil#\hfil\hfil\;\;&\vrule#&
            \;\;\hfil#\hfil\hfil\;\;&\vrule#\cr\tablerule
&&\multispan5\hfil $(\mathfrak{g}, \mathfrak{g}^\tau)$ 
of holomorphic type
 \hfil \ct
&& && $\mathfrak{g}$ && $\mathfrak{g}^\tau$   \ct
&&1&&  $\mathfrak{su}(p,q)$     && $\hphantom{mmi}\mathfrak{s}(\mathfrak{u}(i,j) + \mathfrak{u}(p-i,q-j))$ \ct
&&2&&  $\mathfrak{su}(n,n)$     && $\hphantom{mmi}\mathfrak{so}^*(2n)$ \ct
&&3&&  $\mathfrak{su}(n,n)$     && $\hphantom{mmm}\mathfrak{sp}(n, \mathbb{R})$ \ct
&&4&&  $\mathfrak{so}^*(2n)$    && $\hphantom{mm}\mathfrak{so}^*(2p) + \mathfrak{so}^*(2n-2p)$ \ct
&&5&&  $\mathfrak{so}^*(2n)$    && $\hphantom{mmm}\mathfrak{u}(p,n-p)$ \ct
&&6&&  $\mathfrak{so}(2,n)$     && $\hphantom{m,}\mathfrak{so}(2,p)+\mathfrak{so}(n-p)$ \ct
&&7&&  $\mathfrak{so}(2,2n)$    && $\hphantom{mm}\mathfrak{u}(1,n)$ \ct
&&8&&  $\mathfrak{sp}(n,\mathbb{R})$&& $\hphantom{mmm}\mathfrak{u}(p,n-p)$ \ct
&&9&&  $\mathfrak{sp}(n,\mathbb{R})$&& $\hphantom{mi}\mathfrak{sp}(p,\mathbb{R})+\mathfrak{sp}(n-p,\mathbb{R})$\ct
&&10&&  $\mathfrak{e}_{6(-14)}$  && $\hphantom{,}\mathfrak{so}(10)+\mathfrak{so}(2)$\ct
&&11&&  $\mathfrak{e}_{6(-14)}$  && $\mathfrak{so}^*(10)+\mathfrak{so}(2)$\ct
&&12&&  $\mathfrak{e}_{6(-14)}$  && $\mathfrak{so}(8,2)+\mathfrak{so}(2)$\ct
&&13&&  $\mathfrak{e}_{6(-14)}$  && $\hphantom{i}\mathfrak{su}(5,1)+\mathfrak{sl}(2, \mathbb{R})$\ct
&&14&&  $\mathfrak{e}_{6(-14)}$  && $\mathfrak{su}(4,2)+\mathfrak{su}(2)$\ct
&&15&&  $\mathfrak{e}_{7(-25)}$  && $\hphantom{i}\mathfrak{e}_{6(-78)} +\mathfrak{so}(2)$\ct
&&16&&  $\mathfrak{e}_{7(-25)}$  && $\hphantom{i}\mathfrak{e}_{6(-14)} +\mathfrak{so}(2)$\ct
&&17&&  $\mathfrak{e}_{7(-25)}$  && $\mathfrak{so}(10,2) +\mathfrak{sl}(2,\mathbb{R})$\ct
&&18&&  $\mathfrak{e}_{7(-25)}$  && $\mathfrak{so}^*(12) +\mathfrak{su}(2)$\ct
&&19&&  $\mathfrak{e}_{7(-25)}$  && $\hphantom{mi}\mathfrak{su}(6,2)$\ct
}}
$$
\caption{}
\label{tbl:3.3.1}
\end{table}

\begin{table}[H]
$$
\vbox{
\offinterlineskip
\def\tablerule{\noalign{\hrule}}
\def\ct{&\cr\tablerule}
\halign{\strut#&\vrule#&
            \;\;\hfil#\hfil\hfil\;\;&\vrule#&
            \;\;\hfil#\hfil\hfil\;\;&\vrule#&
            \;\;\hfil#\hfil\hfil\;\;&\vrule#\cr\tablerule
&&\multispan5\hfil $(\mathfrak{g}, \mathfrak{g}^\tau)$ 
of anti-holomorphic type
\hfil \ct
&&&& $\mathfrak{g}$ && \quad $\mathfrak{g}^\tau$   \ct
&&20&&  $\mathfrak{su}(p,q)$    && $\hphantom{mm}\mathfrak{so}(p,q)$ \ct
&&21&&  $\mathfrak{su}(n,n)$     && $\mathfrak{sl}(n,\mathbb{C}) + \mathbb{R}$ \ct
&&22&&  $\mathfrak{su}(2p,2q)$     && \hphantom{mm}$\mathfrak{sp}(p,q)$ \ct
&&23&&  $\mathfrak{so}^*(2n)$    && $\hphantom{mm}\mathfrak{so}(n,\mathbb{C})$ \ct
&&24&&  $\mathfrak{so}^*(4n)$    && $\mathfrak{su}^*(2n) + \mathbb{R}$ \ct
&&25&&  $\mathfrak{so}(2,n)$     && $\hphantom{mm}\,\,\,\mathfrak{so}(1,p)+\mathfrak{so}(1,n-p)$ \ct
&&26&&  $\mathfrak{sp}(n,\mathbb{R})$&& $\hphantom{mm}\,\,\mathfrak{gl}(n,\mathbb{R})$ \ct
&&27&&  $\mathfrak{sp}(2n,\mathbb{R})$&& $\hphantom{mm}\mathfrak{sp}(n,\mathbb{C})$\ct
&&28&&  $\mathfrak{e}_{6(-14)}$  && $\hphantom{mm}\mathfrak{f}_{4(-20)}$\ct
&&29&&  $\mathfrak{e}_{6(-14)}$  && $\hphantom{mm}\mathfrak{sp}(2,2)$\ct
&&30&&  $\mathfrak{e}_{7(-25)}$  && $\phantom{mi}\;\mathfrak{e}_{6(-26)} +\mathfrak{so}(1,1)$\ct
&&31&&  $\mathfrak{e}_{7(-25)}$  && $\hphantom{mmmi}\mathfrak{su}^*(8)$ \ct
}}
$$
\caption{}
\label{tbl:3.3.2}
\end{table}

\subsection{Proof of Theorem~\ref{thm:Hvisible} for $H = K$}
\label{subsec:4.1}
This subsection gives a proof of Theorem~\ref{thm:Hvisible} in the
case where $G$ is a non-compact simple Lie group and $H = K$.
This is an immediate consequence of Lemma~\ref{lem:3.2}
with $\tau = \theta$ (namely, $H = K$) if we show:

\begin{lemma} 
\label{lem:4.1.1}
Suppose $G$ is a non-compact, simply connected, 
 simple Lie group such that $G/K$ is
 a Hermitian symmetric space.
Let $\theta$ be a Cartan involution
 corresponding to $K$.  
Then there exists an involutive automorphism $\sigma$ of $G$
 satisfying the following three conditions:
\begin{align}
& \text{$\sigma$ and $\theta$ commute.}
\label{eqn:2.4.1}
\\
& \text{$\rrank \mathfrak{g}= \rrank \mathfrak{g}^\sigma$.}
\label{eqn:2.4.2}
\\
& \text{$\sigma Z = -Z$.}
\label{eqn:2.4.3}
\end{align}
\end{lemma}

\begin{proof}
The following table gives a choice of 
$\sigma \in \Aut(\mathfrak{g})$
(and hence an automorphism of the simply-connected $G$)
 for each non-compact simple Lie group $G$ of Hermitian type:
%
\begin{table}[H]
$$
\vbox{
\offinterlineskip
\def\tablerule{\noalign{\hrule}}
\halign{\strut#&\vrule#&
            \;\;\hfil#\hfil\hfil\;\;&\vrule#&
            \;\;\hfil#\hfil\hfil\;\;&\vrule#&
                \hfil#\hfil\hfil\,&\vrule#&
            \;\;\hfil#\hfil\hfil\;\;&\vrule#\cr\tablerule
&&\multispan7\hfil $(\mathfrak{g}, \mathfrak{g}^\sigma)$ 
satisfying \eqref{eqn:2.4.1}, \eqref{eqn:2.4.2}
 and \eqref{eqn:2.4.3}
 \hfil &\cr\tablerule
&& ${\mathfrak{g}}$ && $\mathfrak{g}^\sigma$
 &&&& $\rrank \mathfrak{g} = \rrank \mathfrak{g}^\sigma$  &\cr\tablerule
&& ${\mathfrak{su}}(p,q)$ && ${\mathfrak{so}}(p,q)$ &&&& $\min(p,q)$ &\cr\tablerule
&& ${\mathfrak{so}}^*(2n)$ && ${\mathfrak{so}}(n,\mathbb{C})$ &&&& $[\frac{1}{2} n]$ &\cr\tablerule
&& ${\mathfrak{sp}}(n,\mathbb{R})$ && ${\mathfrak{gl}}(n,\mathbb{R})$ &&&& $n$ &\cr\tablerule
&& ${\mathfrak{so}}(2,n)$ && ${\mathfrak{so}}(1,n-1)+{\mathfrak{so}}(1,1)$ 
     &&&& $\min(n,2)$ &\cr\tablerule
&& ${\mathfrak{e}}_{6 (-14)}$  && ${\mathfrak{sp}}(2,2)$&&&& $2$ &\cr\tablerule
&& ${\mathfrak{e}}_{7(-25)}$   && $\mathfrak{su}^*(8)$ &&&& $3$&\cr\tablerule
}}
$$
\caption{}
\label{tbl:4.1.2}
\end{table}
\noindent
All of these pairs $(\mathfrak{g},\mathfrak{g}^\sigma)$ 
appear in Table \ref{tbl:3.3.2},
showing that they are of anti-holomorphic type.
The real rank condition \eqref{eqn:2.4.2} can be verified directly
(see the above Table).
Hence, we have proved Lemma.
\end{proof}

\addtocounter{defin}{1}
\begin{rem}
\label{rem:4.1.3}
The choice of $\sigma$ 
 is not unique.
For example,
 we may choose
 $\mathfrak{g}^\sigma \simeq {\mathfrak{e}}_{6(-26)} + \mathbb{R}$
 in place of $\mathfrak{g}^\sigma \simeq \mathfrak{su}^*(8)$
 if $\mathfrak{g} = {\mathfrak{e}}_{7(-25)}$.
\end{rem}

\subsection{Proof of Theorem~\ref{thm:Gvisible}}
\label{subsec:4.2}
This subsection gives a
proof of Theorem~\ref{thm:Gvisible}
that concerns with the diagonal action of $G$ on 
the direct product
$D \times D$
or $D \times \overline{D}$. 
We shall see that Lemma~\ref{lem:4.1.1} is again a key 
ingredient of the proof
as in the proof of Theorem~\ref{thm:Hvisible} for $H=K$.

Let $G$ be a non-compact, simple Lie group such that $G/K$ is a
Hermitian symmetric space.
We use the letter $\theta'$ in place of the previous $\theta$
to denote the corresponding Cartan
involution of $G$.
Then $\theta(g_1,g_2) := (\theta' g_1, \theta' g_2)$
defines a Cartan involution of the direct product group $G \times G$. 

We define an involutive automorphism $\tau$ of $G \times G$ by
 $\tau(g_1,g_2) := (g_2,g_1)$.
Then $(G \times G)^\tau = \diag (G) $.

\begin{proof}[Proof of Theorem~\ref{thm:Gvisible}]
 \enspace 1) \enspace
Let $\sigma' \in \Aut(G)$ be the involution given in
Lemma~\ref{lem:4.1.1}.
Now, we set
$\sigma(g_1,g_2) := (\sigma' g_1, \sigma' g_2)$.
Obviously, $\tau, \theta$ and $\sigma$ all commute.
Further, $\sigma$ acts on $D \times D$ as an anti-holomorphic
diffeomorphism because so does
$\sigma'$ on $D$
by \eqref{eqn:2.4.3}.
In light that
\begin{alignat*}{3}
& (\mathfrak{g} \oplus \mathfrak{g})^{\tau\theta}
&&   = \{ (X, \theta' X): X \in \mathfrak{g} \} 
&&   \simeq \mathfrak{g},
\\
& (\mathfrak{g} \oplus \mathfrak{g})^{\sigma,\tau\theta}
&&   = \{ (X, \theta' X): X \in \mathfrak{g}^{\sigma'} \} 
&&   \simeq \mathfrak{g}^{\sigma'},
\end{alignat*}
we have 
$\rrank (\mathfrak{g} \oplus \mathfrak{g})^{\tau\theta}
   = \rrank (\mathfrak{g} \oplus \mathfrak{g})^{\sigma,\tau\theta}$
from \eqref{eqn:2.4.2}.
Let $\mathfrak{a}'$ be a maximal abelian subspace of
$\mathfrak{g}^{\sigma',-\theta'}$.
We set
$$
\mathfrak{a} := \{(X,-X): X \in \mathfrak{a}'\}
               = \{(X,\theta'X): X \in \mathfrak{a}'\},
$$
and define a submanifold of $D \times D$ by
$$
 S := \exp \mathfrak{a} \cdot (o,o).
$$
Then Lemma~\ref{lem:3.2} applied to
$(G \times G, D \times D)$ shows that
the diagonal action of $G$ on 
$D \times D = (G \times G)/(K \times K)$
is strongly visible.

2) \enspace
We define an involution $\sigma$ of $G \times G$ by
$\sigma := \tau\theta$, namely,
$\sigma(g_1, g_2) := (\theta' g_2, \theta' g_1)$
 for $g_1, g_2 \in G$.
(This $\sigma$ is different from the one used in 1).)
Then $\sigma$ acts anti-holomorphically on $D \times \overline{D}$.
Obviously,
 $\sigma = \tau\theta$,
 $\tau$ and $\theta$ 
 all commute.
Furthermore,
the rank condition
$
\rrank (\mathfrak{g} \oplus \mathfrak{g})
   ^{\tau\theta}
=
 \rrank (\mathfrak{g} \oplus \mathfrak{g})^{\sigma,\tau\theta}
$
automatically follows from $\sigma = \tau\theta$.
Thus,
it follows from Lemma~\ref{lem:3.2} that
the diagonal action of $G$ on $D \times \overline{D}$ is also
 strongly visible.
Hence, Theorem~\ref{thm:Gvisible} has been proved.
\end{proof}

\section{Visible actions on non-compact $G/K$}
\label{sec:3}
In this section, we give a proof of Theorem~\ref{thm:Hvisible} in the
case where $G$ has no compact factor.
The compact case will be proved in Section~\ref{sec:4}.

\subsection{Existence of anti-holomorphic involutions}
\label{subsec:5.1}
Throughout this section,
let $G$ be a non-compact, simply-connected, simple Lie group of
Hermitian type.
Suppose $\tau$ is an involutive automorphism of $G$ such that
$H = G_0^\tau$.
Then, owing to Lemma~\ref{lem:3.2},
the proof of Theorem~\ref{thm:Hvisible} 
is
reduced to the following:

\begin{lemma}
\label{lem:5.1}
Suppose that 
$\mathfrak{g} = \mathfrak{k} + \mathfrak{p}$ 
is a real simple Lie algebra of Hermitian type
and that $Z$ is a generator of the center 
$c(\mathfrak{k})$ of $\mathfrak{k}$.
Let $\tau$ be an involutive automorphism of\/ $\mathfrak{g}$, 
 commuting with the Cartan involution $\theta$.
Then there exists an involutive automorphism $\sigma$ of $\mathfrak{g}$
 satisfying the following three conditions:
\begin{align}
\label{eqn:3.1.1}
& \text{$\sigma$, $\tau$ and $\theta$ commute with one another.}
\\
\label{eqn:3.1.2}
& \text{$\rrank \mathfrak{g}^{\tau \theta} 
 = \rrank \mathfrak{g}^{\sigma, \tau \theta}$.}
\\
\label{eqn:3.1.3}
& \text{$\sigma Z = -Z$.}
\end{align}
\end{lemma}

The rest of this section will be spent for the proof of 
Lemma~\ref{lem:5.1}. 
We shall divide the proof into the following three cases:

\begin{tabular}{lll}
{Case I.} & $\tau Z = -Z$
          &
(Table~\ref{tbl:3.3.2}, $20 \sim 31$).  \\
{Case II.}& $\tau Z = Z$, $\mathfrak{g}$ is classical
          &
(Table~\ref{tbl:3.3.1}, $1 \sim 9$).    \\
{Case III.}& $\tau Z = Z$, $\mathfrak{g}$ is exceptional\qquad\qquad
           &
(Table~\ref{tbl:3.3.1}, $10 \sim 19$). 
\end{tabular}

\noindent
We have already proved Lemma \ref{lem:5.1}
for $\tau = \theta$ in Subsection \ref{subsec:4.1}
(namely, special cases of Cases II and III).
In the subsequent subsections,
we shall choose $\sigma$ in the following way:

 \begin{tabular}{ll}
{Case I.} 
& Take $\sigma := \tau \theta$ 
 \ (Subsection~\ref{subsec:5.2}). \\
{Case II.} 
& Take $\sigma$ as in Table~\ref{tbl:5.3.1} 
 \ (Subsection~\ref{subsec:5.3}). \\
{Case III.}
& Take $\sigma$ such that $\mathfrak{g}^\sigma \simeq \mathfrak{sp}(2,2)$
if $\mathfrak{g} = \mathfrak{e}_{6(-14)}$, and
 $\mathfrak{g}^\sigma \simeq \mathfrak{su}^*(8)$ \\ 
& if $\mathfrak{g} = \mathfrak{e}_{7(-25)}$
 \     (Subsections~\ref{subsec:5.4}--\ref{subsec:5.6}). 
 \end{tabular}

\subsection{Proof of Lemma~\ref{lem:5.1} 
in Case I}
\label{subsec:5.2}
Suppose $\tau Z = -Z$.
We set $\sigma := \tau \theta$.
Then,
 the conditions
 \eqref{eqn:3.1.1} and \eqref{eqn:3.1.3} are automatically satisfied.
Since
$
\mathfrak{g}^{\sigma, \tau \theta}
=
 \mathfrak{g}^{\tau \theta},
$
 the real rank condition \eqref{eqn:3.1.2}
is obvious.
Thus, Lemma~\ref{lem:5.1} 
in Case I is proved.
\qed

\subsection{Proof of Lemma~\ref{lem:5.1} in Case II}
\label{subsec:5.3}
There are $9$ families of semisimple symmetric pairs
 $(\mathfrak{g}, \mathfrak{g}^\tau)$ in Case II,
namely, the cases $1 \sim 9$ as labeled in Table~\ref{tbl:3.3.1}.
Then, we can take an involutive automorphism
$\sigma$ of $\mathfrak{g}$ as in the following
Table~\ref{tbl:5.3.1}. 
The conditions
 \eqref{eqn:3.1.1} and \eqref{eqn:3.1.3} 
 are clear.
The real rank condition~\eqref{eqn:3.1.2} is verified directly
(see the right column of the table below).
Hence, Lemma~\ref{lem:5.1} in Case II is proved.
\qed

{\small
\begin{table}[H]
$$
\vbox{
\offinterlineskip
\def\tablerule{\noalign{\hrule}}
\def\ct{&\cr\tablerule}
\halign{\strut#&\vrule#&
            \;\hfil#\hfil\hfil\;&\vrule#&
            \,\hfil#\hfil\hfil\,\!&\vrule#&
            \,\hfil#\hfil\hfil\,\!&\vrule#&
            \,\hfil#\hfil\hfil\,\!  &\vrule#\cr\tablerule
&&$(\mathfrak{g},\mathfrak{g}^\tau)$ && $\mathfrak{g}^\sigma$ 
   && $\mathfrak{g}^{\sigma,\tau} = \mathfrak{g}^\sigma \cap \mathfrak{g}^\tau$
&& $\rrank \mathfrak{g}^{\tau\theta}=\rrank\mathfrak{g}^{\sigma,\tau\theta}$
\ct
&&1&&  ${\mathfrak {so}(p,q)}$ 
   && ${\mathfrak{so}(i,j) + \mathfrak{so}(p-i,q-j)}$ 
&& $\min(i,q-j)+\min(p-i,j)$
\ct
&&2&&  ${\mathfrak {so}(n,n)}$    
   && $\mathfrak{so}(n, \mathbb{C})$
&&$n$
\ct
&&3&&  ${\mathfrak {sl}(n,\mathbb{C})+\mathbb{R}}$    
   && $\mathfrak{gl}(n, \mathbb{R})$
&& \raisebox{.1ex}{$[\frac{n}{2}]$}
\ct
&&4&&  ${\mathfrak {so}(n,\mathbb{C})}$
   && ${\mathfrak{so}(p, \mathbb{C}) + \mathfrak{so}(n-p,\mathbb{C})}$
&& $\min(p,n-p)$
\ct
&&5&&  ${\mathfrak {so}(n,\mathbb{C})}$    
   && $\mathfrak{so}(p,n-p)$ 
&& \raisebox{.1ex}{$[\frac{p}{2}] + [\frac{n-p}2] $ }
\ct
&&6&&  ${\mathfrak {so}(1,1)+\mathfrak{so}(1,n-1)}$     
   && ${\mathfrak{so}(1,p)+\mathfrak{so}(n-p-1)}$
&& $\min(2,n-p) $ 
\ct
&&7&&  ${\mathfrak {so}(1,n)+\mathfrak{so}(1,n)}$     
   && $\mathfrak{so}(1,n)$ 
&& $1 $ 
\ct
&&8&&  ${\mathfrak {gl}(n,\mathbb{R})}$
   && $\mathfrak{so}(p,n-p)$
&& $n $
\ct
&&9&&  ${\mathfrak {gl}(n,\mathbb{R})}$
   && ${\mathfrak{gl}(p,\mathbb{R})+\mathfrak{gl}(n-p,\mathbb{R})}$
&&$\min(p,n-p) $
\ct
}}
$$
\caption{}
\label{tbl:5.3.1}
\end{table}
}

\subsection{$\epsilon$-family of symmetric pairs}
\label{subsec:5.4}
We shall prove Lemma~\ref{lem:5.1} 
in Case III.
We shall take $\sigma$ so that
 $\mathfrak{g}^\sigma \simeq \mathfrak{sp}(2,2)$ for
 $\mathfrak{g} = \mathfrak{e}_{6(-14)}$,
 and
 $\mathfrak{g}^\sigma \simeq \mathfrak{su}^*(8)$
 for $\mathfrak{g} = \mathfrak{e}_{7(-25)}$.
The non-trivial part is to prove that we can take $\sigma$
 such that $\sigma \tau = \tau \sigma$.
(Two involutions do not always commute.
 See Subsection~\ref{subsec:7.6} for  counterexamples in classical cases.)

We have already proved that Theorem~\ref{thm:Hvisible} holds if 
$\tau = \theta$ (see Subsection~\ref{subsec:4.1})
or if $\tau$ is of anti-holomorphic type (see Subsection~\ref{subsec:5.2}).
Building on these cases,
we shall give a proof of the remaining cases,
that is, Lemma~\ref{lem:5.1} in Case III.
For this,
we set up to make new pairs from old. First,
 we recall quickly the notion of $\epsilon$-families
 of symmetric pairs \cite{xosadv},
which enables us to
 avoid tedious computations for exceptional groups.
Our approach below might be of some use for
 a systematic study of three involutions $(\sigma, \tau, \theta)$
 of complex simple Lie algebras (cf. \cite{xmahgl}).

Let ${\mathfrak{g}}$ be a semisimple Lie algebra, 
 $\tau$  an involutive automorphism of ${\mathfrak{g}}$,
 and $\theta$  a Cartan involution of ${\mathfrak{g}}$ commuting with $\tau$.
Fix a maximal abelian subspace $\mathfrak{a}$
 of $\mathfrak{g}^{-\theta, -\tau}$.
For $\lambda \in \mathfrak{a}^*$,
we define
$\mathfrak{g}(\mathfrak{a};\lambda) := \{ X \in \mathfrak{g}:
  \operatorname{ad}(H)X = \lambda(H)X \text{ for }
  H \in \mathfrak{a} \}
$,
and set
$
\Sigma(\mathfrak{a}) \equiv \Sigma(\mathfrak{g},\mathfrak{a})
  := \{ \lambda \in \mathfrak{a}^* \backslash \{0\} : 
     \mathfrak{g}  (\mathfrak{a};\lambda) \ne \{0\} \}
$.
Rossmann proved that
 $\Sigma(\mathfrak{a}) $ satisfies the axiom 
 of root system (\cite[Theorem 5]{xros}).   
We say a map 
$$
\epsilon : \Sigma(\mathfrak{a}) \cup \{0\} \to \{ 1, -1\}
$$
 is a \textit{signature of\/ $\Sigma(\mathfrak{a})$} 
 if $\epsilon$ satisfies
 $\epsilon(\alpha + \beta) = \epsilon(\alpha) \epsilon(\beta)$ for any
 $\alpha,\, \beta,\, \alpha + \beta \in \Sigma(\mathfrak{a})$,
 $\epsilon(-\alpha) = \epsilon(\alpha)$ for any $\alpha \in \Sigma(\mathfrak{a})$
 and $\epsilon(0) = 1$. 
To a signature $\epsilon$ of $\Sigma(\mathfrak{a})$, 
 we associate an involution $\tau_\epsilon$ of $\mathfrak{g}$ defined
 by 
$$
\tau_\epsilon (X) := \epsilon(\lambda) \tau(X)
\text{\quad for \ }
X \in {\mathfrak{g}}(\mathfrak{a};\lambda), \ 
 \lambda \in \Sigma(\mathfrak{a}) \cup \{ 0\}.
$$
Then $\tau_\epsilon$ defines another symmetric pair 
 $({\mathfrak{g}},{\mathfrak{h}}_\epsilon)$.
The set 
$$
    F\left(({\mathfrak{g}},{\mathfrak{h}})\right)
    :=
    \set{({\mathfrak{g}},{\mathfrak{h}}_\epsilon)}{\epsilon \text{ is a signature of }\Sigma(\mathfrak{a})}
$$
 is said to be an {\it $\epsilon$-family} 
 {\it {of symmetric pairs}} (\cite[\S 6]{xosadv}).
This set is also referred to as
 $K_\epsilon$-\textit{family of symmetric pairs}
 if $\tau = \theta$.
For example, 
 $\set{(\mathfrak{sl}(n, \mathbb{R}), \mathfrak{so}(p,n-p))}{0 \le p \le n}$
 forms a $K_\epsilon$-family.

To make new pairs $(\sigma, \tau)$ from old,
 the following result is useful: 
\begin{lemma}
\label{lem:5.4}
Let $\tau$ and $\sigma$ be involutive automorphisms of\/ $\mathfrak{g}$.
\itm{1)}
If the pair $(\sigma, \tau)$ satisfies \eqref{eqn:3.1.1}
 and \eqref{eqn:3.1.2},
then so does $(\sigma, \tau_\epsilon)$ for any $\tau_\epsilon$
 (the choice of\/ $\mathfrak{a}$ is specified in the proof).
\itm{2)}
If the pair $(\sigma, \tau)$ satisfies \eqref{eqn:3.1.1},
then so does $(\sigma, \tau \theta)$.
\end{lemma}
\begin{proof}
1) \enspace
Take a maximal abelian subspace $\mathfrak{a}$ in
$$
   \mathfrak{g}^{-\theta, \sigma, \tau \theta}
   =
   \mathfrak{g}^{-\theta, \sigma, -\tau}
   =
   \set{Y \in \mathfrak{g}}{(-\theta) Y = \sigma Y = (-\tau) Y = Y} 
    .
$$
Then, $\mathfrak{a}$
 is also a maximal abelian subspace of $\mathfrak{g}^{-\theta, -\tau}$
 by the rank condition \eqref{eqn:3.1.2}.
Let $\tau_\epsilon$ be an involutive automorphism of $\mathfrak{g}$ associated to 
 a signature $\epsilon$ of $\Sigma(\mathfrak{g}, \mathfrak{a})$.

To see $\sigma \tau_\epsilon = \tau_\epsilon \sigma$,
 we take an arbitrary root vector
$X \in {\mathfrak{g}}(\mathfrak{a};\lambda)$.
Then we have
$$
   \sigma \tau_\epsilon(X) = \sigma(\epsilon(\lambda) \tau(X))
                              =  \epsilon(\lambda) \sigma(\tau(X))  .
$$
On the other hand,
 $\sigma X \in {\mathfrak{g}}(\mathfrak{a};\sigma \lambda)
 ={\mathfrak{g}}(\mathfrak{a};\lambda)$
 because $\sigma|_{\mathfrak{a}} = \mathrm{id}$.
Thus,
$$
   \tau_\epsilon \sigma (X)  =  \epsilon(\lambda) \tau(\sigma (X))
                               =  \epsilon(\lambda) \sigma(\tau(X)) 
                                   .
$$
This proves $\sigma \tau_\epsilon = \tau_\epsilon \sigma$ on 
 $\mathfrak{g}(\mathfrak{a};\lambda)$ for all $\lambda$.
Hence,
 the pair $(\sigma, \tau_\epsilon)$ satisfies \eqref{eqn:3.1.1}.

As $\epsilon(0) = 1$,
 we have $\tau_\epsilon(Y) = \tau(Y) = -Y$ if
 $Y \in \mathfrak{a} \ (\subset \mathfrak{g}(\mathfrak{a};0))$.
This implies 
$$
   \mathfrak{a} \subset
   \mathfrak{g}^{-\theta, \sigma, \tau_\epsilon \theta}
   =
   \mathfrak{g}^{-\theta, \sigma, -\tau_\epsilon}
   =
   \set{Y \in \mathfrak{g}}{(-\theta) Y = \sigma Y = (-\tau_\epsilon) Y = Y}
    .
$$
Since $\rrank \mathfrak{g}^{\sigma, \tau_\epsilon \theta}$ 
 is  by definition the dimension of a maximal abelian subspace
 contained in $\mathfrak{g}^{-\theta, \sigma, \tau_\epsilon \theta}$,
 we have 
$$
\rrank \mathfrak{g}^{\sigma, \tau_\epsilon \theta} 
\ge \dim \mathfrak{a} = \rrank \mathfrak{g}^{\tau \theta} 
 = \rrank \mathfrak{g}^{\tau_\epsilon \theta} 
 \ge \rrank \mathfrak{g}^{\sigma, \tau_\epsilon \theta}  . 
$$
Therefore,
 the pair $(\sigma, \tau_\epsilon)$ satisfies \eqref{eqn:3.1.2}.
\itm{2)}
Obvious from the definition.
\end{proof}

\subsection{New pairs $(\sigma,\tau)$ from old}
\label{subsec:5.5}
We have already proved Lemma~\ref{lem:5.1} 
in the following two cases:
\itm{i)}
 $\tau = \theta$ (equivalent to Lemma~\ref{lem:4.1.1}).
\itm{ii)}
 $\tau$ satisfies $\tau Z = -Z$ 
(Subsection~\ref{subsec:5.2}).
\newline
The lemma below gives a coherent understanding of the set of
involutions for which Lemma \ref{lem:5.1} holds.
\begin{lemma}
\label{lem:5.5}
Let $\tau$ be an involutive automorphism of\/ $\mathfrak{g}$ commuting
 with a Cartan involution $\theta$.
\itm{1)}
If Lemma~\ref{lem:5.1} 
holds for $\tau$,
 then so does it for $\tau_\epsilon$
associated to any signature $\epsilon$ of
 $\Sigma(\mathfrak{g},\mathfrak{a})$. 
To be more precise,
let $\sigma$ be an automorphism of $\mathfrak{g}$ satisfying
\eqref{eqn:3.1.1}, \eqref{eqn:3.1.2} and \eqref{eqn:3.1.3}, and
 $\mathfrak{a}$ be a maximal abelian subspace in
$\mathfrak{g}^{-\theta,\sigma,-\tau}$.
Then, $\tau_\epsilon\theta = \theta\tau_\epsilon$ and the same
 $\sigma$ satisfies \eqref{eqn:3.1.1}, \eqref{eqn:3.1.2} and
 \eqref{eqn:3.1.3} for $\tau_\epsilon$.
\itm{2)}
If Lemma~\ref{lem:5.1} 
holds for $\tau$ by taking an involution $\sigma$,
 then so does it for $\tau \theta$,
 provided $\rrank \mathfrak{g}^\tau = \rrank \mathfrak{g}^{\sigma, \tau}$.
Namely, the same $\sigma$ satisfies \eqref{eqn:3.1.1},
 \eqref{eqn:3.1.2} and \eqref{eqn:3.1.3} for $\tau\theta$.
\end{lemma}

\begin{proof}
Readily follows from
Lemma~\ref{lem:5.4}. 
\end{proof}

The next subsection shows that
 Lemma~\ref{lem:5.1} 
in  the exceptional cases 
(namely, (10) $\sim$ (19) in Table~\ref{tbl:3.3.1})
 is reduced to (i) or (ii) by an iterating application of
 Lemma~\ref{lem:5.5}. 

\subsection{Proof of Lemma~\ref{lem:5.1} in Case~III}
\label{subsec:5.6}
In terms of the labels of symmetric pairs in 
Tables~\ref{tbl:3.3.1} and \ref{tbl:3.3.2},
the scheme of the proof is described in the following diagrams:
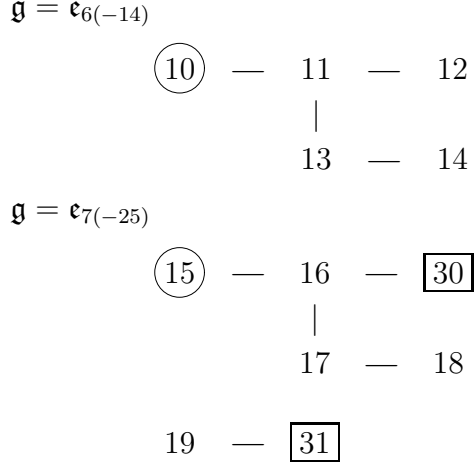
\begin{alignat*}{2}
&\mathfrak{g} = \mathfrak{e}_{6(-14)}
\\[1ex]
&&&\begin{array}{ccccc}
\put(5.8,4){\circle{19}}10 \ 
     & \rotatebox[origin=c]{90}{$|$} & \ 11 \ & \rotatebox[origin=c]{90}{$|$} & \ 12 \ 
    \\[.5ex]
                         &&|
    \\[.5ex]
                         && 13 & \rotatebox[origin=c]{90}{$|$} & 14
 \end{array}
\\
&\mathfrak{g} = \mathfrak{e}_{7(-25)}
\\[1ex]
&&&\begin{array}{ccccc}
\put(5.8,4){\circle{19}}15 \ 
        & \rotatebox[origin=c]{90}{$|$} & 16 & \rotatebox[origin=c]{90}{$|$} 
        &  \fbox{30}
    \\[.5ex]
                         &&|
    \\[.5ex]
                         &&17 & \rotatebox[origin=c]{90}{$|$} & 18
    \\              
    \\
          19 \     & \rotatebox[origin=c]{90}{$|$} 
                 & \fbox{31}
  \end{array}
\hspace*{10em}
\end{alignat*}
Here, the symmetric pairs connected by
horizontal path mean that they belong to
the same $\epsilon$-family.
That is,
we recall from \cite{xosadv} that
the following is a list of
 an $\epsilon$-family  of symmetric pairs:
\itm{$\mathfrak{g} =\mathfrak{e}_{6(-14)}$}:
$\{(10), (11), (12)\}, \
\{(13), (14)\}$. 
\itm{$\mathfrak{g} = \mathfrak{e}_{7(-25)}$}:
$\{(15), (16), (30)\}, \
\{(17), (18)\}, \
\{(19), (31)\}$.

The circle
(i.e.\ (10) and (15)) means that $\tau = \theta$,
while the  box
(i.e.\ (30) and (31)) means that $\tau$ is of anti-holomorphic type. 

Since Lemma~\ref{lem:5.1} 
holds for (10), (15) ($\tau = \theta$ case)
 and also for (30), (31) ($\tau Z = -Z$ case),
 so does it for any member of
$\{(10), (11), (12)\}, \
\{(15), (16)$, $(30)\}, \
\{(19), (31)\}$
by Lemma~\ref{lem:5.5}~(1).  

A next step is an observation:
\itm{$\mathfrak{e}_{6(-14)}$}:
Lemma~\ref{lem:5.1} 
holds for (11) by taking $\mathfrak{g}^\sigma \simeq \mathfrak{sp}(2,2)$.
If $(\mathfrak{g}, \mathfrak{g}^\tau) \simeq$
 (11),
\newline\indent\indent
 then $(\mathfrak{g}, \mathfrak{g}^{\tau \theta}) \simeq$ (13) and
 $(\mathfrak{g}^\tau, \mathfrak{g}^{\sigma, \tau})
 \simeq (\mathfrak{so}^*(10)+\mathfrak{so}(2), \mathfrak{sp}(2, \mathbb{C}))$.
\itm{$\mathfrak{e}_{7(-25)}$}:
Lemma~\ref{lem:5.1} 
holds for (16) by taking
 $\mathfrak{g}^\sigma \simeq \mathfrak{su}^*(8)$.
If $(\mathfrak{g}, \mathfrak{g}^\tau) \simeq$
 (16),
\newline\indent\indent
 then $(\mathfrak{g}, \mathfrak{g}^{\tau \theta}) \simeq$  (17) and
 $(\mathfrak{g}^\tau, \mathfrak{g}^{\sigma, \tau})
 \simeq (\mathfrak{e}_{6(-14)}+\mathfrak{so}(2), \mathfrak{sp}(2,2))$.

Here, the proof of the above isomorphisms concerning
 $(\mathfrak{g}^\tau, \mathfrak{g}^{\sigma, \tau})$ 
is straightforward because we know $\sigma \tau = \tau \sigma$.

Then, Lemma~\ref{lem:5.1} 
holds for (13) and (17) by Lemma~\ref{lem:5.5}~(2).  
In turn,
 Lemma~\ref{lem:5.1} 
holds for  any member of $\{(13), (14)\}$ and $\{(17), (18)\}$
 by using again Lemma~\ref{lem:5.5}~(1).  
 
This proves Lemma~\ref{lem:5.1}  
 for all exceptional cases (10) $\sim$ (19), namely, in Case III.
Thus, we have finished the proof of Lemma~\ref{lem:5.1}. 
\qed

\section{Visible actions on compact $G/K$}
\label{sec:4}
This section gives a proof of Theorem~\ref{thm:Hvisible} in the case
where $D$ is a compact symmetric space.
This is reduced to the following two cases:

Case~I. \quad
$D = G_U / K$.

Case~II.\quad
$D = (G_U \times G_U) / (K_1 \times K_2), \
 H_U = \diag (G_U)$.

\noindent
Here, $G_U$ is a connected compact simple Lie group,
and $G_U / K, G_U / K_1$ and $G_U / K_2$ are compact irreducible
Hermitian symmetric spaces.
(Instead of the letters $G$ and $H$,
we shall use $G_U$ and $H_U$ to emphasize compactness in 
Section~\ref{sec:4}). 

Theorem~\ref{thm:Hvisible} in Cases~I and II will be proved in
Subsections~\ref{subsec:cpt1} and 
\ref{subsec:cpt2}, respectively.
Together with the non-compact case proved in Section~\ref{sec:3},
the proof of Theorem \ref{thm:Hvisible} will be completed.

\subsection{Existence of anti-holomorphic involutions}
\label{subsec:4.1new}
Without loss of generality,
we may and do assume that $G_U$ is simply connected.
We denote by $\tau$ and $\theta$ the involutive automorphisms of $G_U$
such that $(G_U)^\tau = H_U$ and $(G_U)^\theta = K$,
respectively. 
Since $G_U$ is simply connected,
both $H_U$ and $K$ are automatically  connected.

For $g \in G_U$, we define an involution $\tau^g$ by
$$
   \tau^g(x) := g \tau(g^{-1} x g) g^{-1} \quad (x \in G_U).
$$
Here is a key lemma:

\begin{lemma}
\label{lem:7.5}
Suppose we are in the above setting.
Then,
there exists an involutive automorphism $\sigma$ of $G_U$
 satisfying the following three conditions
 (by an abuse of notation,
 we replace $\tau^g$ with $\tau$ below
by taking some $g \in G_U$ if necessary):
\begin{align}
& \text{$\sigma \theta = \theta \sigma$, $\sigma \tau = \tau \sigma$.}
\label{eqn:7.5.1}
\\
& \text{$(\mathfrak{g}_U)^{\sigma, -\tau, -\theta}$
 contains a maximal abelian subspace in  $(\mathfrak{g}_U)^{-\tau, -\theta}$.}
\label{eqn:7.5.3}
\\
& \text{The induced action of $\sigma$ on $D = G_U/K$ is anti-holomorphic.}
\label{eqn:7.5.2}
\end{align}
\end{lemma}

\begin{remark}
Lemma~\ref{lem:7.5} is a compact case counterpart of
Lemma~\ref{lem:5.1}.
In contrast to the condition~\eqref{eqn:3.1.1}
in the non-compact case,
we have not required $\tau \theta = \theta \tau$ here.
In fact,
different from the non-compact case,
it may happen that 
$\tau^g \theta \ne \theta \tau^g$
for all $g \in G_U$
(see Type~II in Subsection~\ref{subsec:7.6}).
Nevertheless, Lemma \ref{lem:7.5} asserts that one can find $\sigma$
that commutes with $\theta$ and $\tau$ simultaneously.
\end{remark}

We shall divide the proof into the following cases:
\newline
Type I:\ \ $\tau^g \theta = \theta \tau^g$ for some $g \in G_U$.
\newline
Type II:\ $\tau^g \theta \neq \theta \tau^g$ for any $g \in G_U$.

\subsection{Proof of Lemma~\ref{lem:7.5} in Type I} 
\label{subsec:7.6}
This subsection gives a proof of Lemma~\ref{lem:7.5} in Type I.
Type I parallels the corresponding result
 (see Lemma~\ref{lem:5.1})
 for the non-compact Riemannian symmetric pair $(G, K)$ dual to $(G_U, K)$.

Let $G_\mathbb{C}$ be a complexification of $G_U$.
Since $G_U$ is simply connected,
$G_\mathbb{C}$ is also simply connected.
Therefore, any automorphism of $G_U$ extends to a holomorphic
automorphism of the complex Lie group $G_\mathbb{C}$.
For $\tau, \theta, \ldots \in \Aut (G_U)$,
we shall use the same letters $\tau, \theta, \ldots$ to denote the
holomorphic extensions 
$\in \Aut (G_{\mathbb{C}})$,
and also the differentials
$\in \Aut (\mathfrak{g}_{\mathbb{C}})$.

Let $G$ be a connected subgroup of $G_\mathbb{C}$ with Lie algebra
\begin{equation}
\mathfrak{g}:= \mathfrak{k} + \sqrt{-1} (\mathfrak{g}_U)^{-\theta}
 = (\mathfrak{g}_U)^\theta + \sqrt{-1} (\mathfrak{g}_U)^{-\theta}  .
\label{eqn:7.6.1}
\end{equation}
Then $G$ is a non-compact simple Lie group with
 a maximal compact subgroup $K$.

In Type I, we shall simply write $\tau$ for $\tau^g$
and may assume $\tau \theta = \theta \tau$.
Then, the decomposition \eqref{eqn:7.6.1} 
is invariant by $\tau$,
and consequently, the holomorphic extension
 $\tau \in \Aut(G_\mathbb{C})$ stabilizes $G$.
Now, we take $\sigma \in \Aut(\mathfrak{g})$ as in Lemma~\ref{lem:5.1},
 extend it holomorphically on $G_\mathbb{C}$,
 and restrict it to $G_U$  (we use the same letter $\sigma$). 

Then, $\sigma\theta = \theta\sigma$ and $\sigma\tau = \tau\sigma$ hold
on $G_U$ because so do they on $G$.
It follows from $\sigma\theta = \theta\sigma$ that $\sigma$ induces a
 diffeomorphism of $D = G_U/K$.
Furthermore,
this is anti-holomorphic, because $\sigma Z = -Z$ and the complex
 structure on $G_U/K$ is given by the left $G_U$-translation of 
$\operatorname{Ad} (\exp (\frac{\pi}{2} Z))$
as is the case of the non-compact Hermitian symmetric space $G/K$
(see Subsection~\ref{subsec:3.3}).

Since the condition \eqref{eqn:3.1.2} means that
$\mathfrak{g}^{\sigma,-\tau,-\theta}$ 
contains a maximal abelian subspace of 
$\mathfrak{g}^{-\tau,-\theta}$, 
the condition \eqref{eqn:7.5.3} follows from 
$$
   (\mathfrak{g}_U)^{-\tau,-\theta}
  =
   \sqrt{-1}\mathfrak{g}^{-\tau,-\theta} 
  \quad\text{and}\quad
  (\mathfrak{g}_U)^{\sigma,-\tau,-\theta}
  = 
   \sqrt{-1} \mathfrak{g}^{\sigma,-\tau,-\theta}
 .
$$
Therefore, all the conditions \eqref{eqn:7.5.1} -- \eqref{eqn:7.5.2} are
verified.
Thus, we have proved Lemma~\ref{lem:7.5} 
for Type~I.
\qed

\subsection{Proof of Lemma \ref{lem:7.5} in Type II}
\label{subsec:41TypeII}

In contrast to Type I,
  $\tau$ (or any of its conjugation)
 cannot stabilize a non-compact real form $G$ of $G_{\mathbb{C}}$ in Type~II.
Thus, we cannot reduce Type~II to the non-compact
results in Section~\ref{sec:3}.
Instead, our strategy here is to find a \lq\lq large\rq\rq\
 subalgebra, say $\mathfrak{g}_U'$, of $\mathfrak{g}_U$,
 such that $\tau$ commutes with $\theta$ when restricted to $\mathfrak{g}_U'$.
The definition of $\mathfrak{g}_U'$ and a
precise formulation of ``largeness'' will be given in
Claim~\ref{claim:7.6}. 

There are two cases up to conjugation for Type~II:
\newline  Type II-1 \
$(\mathfrak{g}_U, \mathfrak{g}_U^\tau, \mathfrak{g}_U^\theta) =
 (\mathfrak{su}(2n), \mathfrak{sp}(n),
 \mathfrak{s u}(2p'+1) + \mathfrak{s u}(2q'+1) + \sqrt{-1}\mathbb{R})$,
\newline Type II-2 \
$(\mathfrak{g}_U, \mathfrak{g}_U^\tau, \mathfrak{g}_U^\theta) =
 (\mathfrak{so}(2n), \mathfrak{so}(2p'+1) + \mathfrak{so}(2q'+1), \mathfrak{u}(n))$,
\newline
 where $n = p'+q'+1$.

The proof for Type~II-2 follows from 
Type II-1 by switching the role of $\tau$ and $\theta$ in 
$\mathfrak{so}(2n) = \mathfrak{su}(2n) \cap \mathfrak{gl}(2n, \mathbb{R})$.
Therefore, we shall consider mostly Type II-1,
but for the convenience of the reader, we sometimes supply the formula for 
Type~II-2 in addition.

By using matrix realization,
we suppose $\mathfrak{g}_U^\theta$ is a subalgebra of $\mathfrak{g}_U$
defined by
\newline
$$
\theta(W): = g_{0}\ W \ g_{0}^{-1},
$$
where
$$g_{0}: = \diag (\underbrace{1,\dots ,1}_{p'},
\underbrace{-1,\dots , -1}_{q'+1}, 
\underbrace{1, \dots, 1}_{p'+1},
\underbrace{-1, \dots ,-1}_{q'})
\ \ \in \ \ GL(2n,\mathbb{R}).
$$
Furthermore, by taking conjugation of $\tau$ by $G_U$ if necessary,
we may and do assume that $\mathfrak{g}_U^\tau$ is a subalgebra of 
$\mathfrak{g}_U$ defined by
$$
X \mapsto J_n \overline{X} J_n^{-1},
$$
where we set
$$
J_n
:= \begin{pmatrix}  0 & -I_{n}\\
                   I_{n} & 0
  \end{pmatrix}
\in GL(2n,\mathbb{R}).
$$
Then, we have
\begin{equation}
    \mathfrak{g}_U^{\tau, \theta} \equiv
    \mathfrak{g}_U^{\tau} \cap \mathfrak{g}_U^{\theta} \simeq 
         \begin{cases}
          \mathfrak{sp}(p') + \mathfrak{sp}(q') + \mathfrak{u}(1)
           & \text{(Type II-1)},
           \\
          \mathfrak{u}(p') + \mathfrak{u}(q')
           & \text{(Type II-2)}.
          \end{cases}
\label{eqn:7.6.6}
\end{equation}

The non-commutativity $\tau\theta \ne \theta\tau$ arises from the odd
parity of $2p'+1$ and $2q'+1$ in the both cases,
and is reflected by the
fact that neither 
$(\mathfrak{g}_U^\tau,\mathfrak{g}_U^{\tau,\theta})$
nor
$(\mathfrak{g}_U^\theta,\mathfrak{g}_U^{\tau,\theta})$
is a symmetric pair.
The idea of the following claim is to pull $2p'$ and $2q'$ out of
$2p'+1$ and $2q'+1$.

\begin{claimsubsec}
\label{claim:7.6}
We can realize the Lie algebra:
$$
\mathfrak{g}_U' := \begin{cases}
                   \mathfrak{su}(2n-2) & \text{\rm{(Type II-1)},}
                   \\
                   \mathfrak{so}(2n-2) & \text{\rm{(Type II-2)},}
             \end{cases}
$$
as a subalgebra of $\mathfrak{g}_U$ such that the following three
 conditions are satisfied:
\begin{align}
& \text{Both $\tau$ and $\theta$ stabilize $\mathfrak{g}_U'$.}
\label{eqn:7.6.3}
\\
& \text{$\tau|_{\mathfrak{g}_U'}$ and $\theta|_{\mathfrak{g}_U'}$ commute.}
\label{eqn:7.6.4}
\\
& \text{$(\mathfrak{g}_U)^{-\tau, -\theta} = (\mathfrak{g}_U')^{-\tau, -\theta}$.}
\label{eqn:7.6.5}
\end{align}
\end{claimsubsec}

\begin{proof}
We consider the subspace
$\mathbb{C}^{2n-2} = \mathbb{C}^{2(p'+q')} \subset \mathbb{C}^{2n}$
corresponding to the partition
$$
2n = p'+1+q'+p'+1+q',
$$
and embed $\mathfrak{g}_U'$ in $\mathfrak{g}_U$ 
accordingly.
Then, clearly \eqref{eqn:7.6.3} and \eqref{eqn:7.6.4} hold,
and the triple
$(\mathfrak{g}_U', (\mathfrak{g}_U')^\tau, (\mathfrak{g}_U')^\theta)$
of Lie algebras is given by
\begin{alignat*}{2}
&  (\mathfrak{su}(2n-2), \mathfrak{sp}(n-1), 
   \mathfrak{su}(2p') + \mathfrak{su}(2q') + \sqrt{-1} \mathbb{R})
&&\quad \text{(Type II-1)},
\\
&  (\mathfrak{so}(2n-2), \mathfrak{so}(2p') + \mathfrak{so}(2q'),
   \mathfrak{u}(n-1))
&&\quad \text{(Type II-2)}.
\end{alignat*}
We recall $n = p' + q' +1$.
For $A \in M(p',q';\mathbb{C})$, we set
$$
 X(A): = \begin{pmatrix} 0 & 0 & A \\
                         0 & 0 & 0 \\
                    -A^{*} & 0 & 0 
         \end{pmatrix},
 \
 Y(A): = \begin{pmatrix} 0 & 0 & A \\ 
                         0 & 0 & 0 \\
                    -^{t}\!A & 0 & 0    
         \end{pmatrix}
 \in M(n,\mathbb{C}).
$$
Then, by a simple matrix computation we have
$$
 (\mathfrak{g}_{U})^{-\tau, -\theta} =
 \begin{cases} 
     \biggl\{ \begin{pmatrix} X(A) & Y(B) \\
                 \overline{Y(B)} & -\overline{X(A)}
            \end{pmatrix} & : A, B \in M(p',q' ; \mathbb{C})
     \biggr\} \  \text{(Type II-1)}, \\[3ex]
     \biggl\{ \begin{pmatrix} X(A) & Y(B) \\
                          -Y(B) & -X(A)
            \end{pmatrix} \kern-2em& : A, B \in M(p',q' ; \mathbb{R})
     \biggr\} \  \text{(Type~II-2)}.
 \end{cases}
$$
Now, \eqref{eqn:7.6.5} is clear.
Thus, Claim \ref{claim:7.6} is proved.
\end{proof}

Let us return to the proof of Lemma~\ref{lem:7.5} 
in Type II.

We define
$\sigma \in \Aut(\mathfrak{g}_U)$ by
$$
\sigma(W): = I_{n,n} \ \overline{W} \ I_{n,n}^{-1},
$$
where
$$
I_{n,n} : = \begin{pmatrix} I_{n} & 0 \\
                            0     & -I_{n}
            \end{pmatrix} \in \ GL(2n,\mathbb{R}).
$$
Then we have
 $\sigma\theta = \theta\sigma$, $\sigma\tau = \tau\sigma$ and 
\begin{equation}
 \mathfrak{g}_U^\sigma \simeq \begin{cases}
                           \mathfrak{so}(2n) & \quad \text{(Type II-1)},
                           \\
                           \mathfrak{so}(n)+\mathfrak{so}(n) &\quad \text{(Type II-2)},
                         \end{cases}
\label{eqn:7.6.2}
\end{equation}
$$
 (\mathfrak{g}_{U})^{\sigma,-\tau, -\theta} =
 \begin{cases} 
     \biggl\{\begin{pmatrix} X(A) & Y(B) \\
                 \overline{Y(B)} & -\overline{X(A)}
             \end{pmatrix} \kern-.7em& : A, \sqrt{-1}B \in M(p',q' ; \mathbb{R})
     \biggr\} \  \text{(Type II-1)}, \\[3ex]
     \biggl\{\begin{pmatrix} X(A) & 0 \\
                            0 & -X(A)
             \end{pmatrix} \kern-.7em& : A \in M(p',q' ; \mathbb{R})
     \biggr\}
       \text{\phantom{$,\sqrt{-1}B$} }  
            \text{(Type~II-2)}.
 \end{cases}
$$
We note that $\sigma$ acts anti-holomorphically on both
$G_U/K$ and its complex submanifold
$G'_U/K' := G'_U / G'_U \cap K$
(see Table \ref{tbl:3.3.2}).

Now, we consider
 $(\mathfrak{g}'_U, \theta|_{\mathfrak{g}_U'}, \tau|_{\mathfrak{g}_U'})$.
This is of Type I by \eqref{eqn:7.6.4}, 
for which Lemma~\ref{lem:7.5} has been already proved in 
Subsection \ref{subsec:7.6}. 
In fact, the restriction
$\sigma|_{\mathfrak{g}'_U}$ of the above choice of $\sigma$
satisfies the conclusion of 
Lemma~\ref{lem:7.5}.
In particular,
if we take a maximal abelian subspace $\mathfrak{t}$ of
$(\mathfrak{g}'_U)^{\sigma,-\tau,-\theta}$,
then it is also a maximal abelian subspace of
$(\mathfrak{g}'_U)^{-\tau,-\theta}$.

Next, we consider $(\mathfrak{g},\theta,\tau)$ which is of Type~II.
In view of \eqref{eqn:7.6.5}, 
$\mathfrak{t}$ is also a maximal abelian subspace of
$(\mathfrak{g}_U)^{-\tau,-\theta}$.
Hence, the condition~\eqref{eqn:7.5.3} is satisfied.
We have already seen \eqref{eqn:7.5.1} and \eqref{eqn:7.5.2}.
Thus, we have proved Lemma~\ref{lem:7.5} 
in Type II.
\qed

\subsection{Stability of $H_U$-orbits}
\label{subsec:cpt1}
We are ready to complete the proof of Theorem~\ref{thm:Hvisible} 
in Case I, namely, for a
compact simple $G_U$.
Here is a compact case counterpart of
Lemma~\ref{lem:3.2}:

\begin{lemm}
\label{lem:cptvisible}
Suppose that three involutive automorphisms,
$\tau, \theta$ and $\sigma$ of $G_U$ satisfy the conditions
\eqref{eqn:7.5.1}, \eqref{eqn:7.5.3} and \eqref{eqn:7.5.2}.
Then,\\
{\upshape 1)}\enspace
For any $x \in D = G_U/K$,
there exists $g \in H_U$ such that $\sigma(x) = g \cdot x$. 
In particular,
each $H_U$-orbit on $D$ is preserved by $\sigma$.
\\
{\upshape 2)}\enspace
Take a maximal abelian subspace $\mathfrak{t}$ in
$(\mathfrak{g}_U)^{\sigma,-\tau,-\theta}$,
and we define a submanifold $S$ of $D = G_U / K$
by $S := (\exp \mathfrak{t})K$.
Then, $S$ meets every $H_U$-orbit in $D$,
and $\sigma |_S = \mathrm{id}$.
\end{lemm}

The proof parallels that of Lemma \ref{lem:3.2}.
For this, all we need now is the following lemma, which is
 a compact analog of a generalized Cartan decomposition $G= H A K$
 (see \eqref{eqn:3.2.3}).
\begin{lemm}
\label{lem:7.4}
Let $G_U$ be a semisimple connected compact Lie group
 with Lie algebra $\mathfrak{g}_U$,
 and $\tau$ and $\theta$ two involutive automorphisms of $G_U$.
Take a maximal abelian subspace $\mathfrak{t}$ in $(\mathfrak{g}_U)^{-\tau, -\theta}$,
 and let $T$ be the analytic subgroup of $G_U$ with Lie algebra $\mathfrak{t}$.
Then we have
$$
 G_U = (G_U)^\tau_0 \ T\ (G_U)^\theta_0 \; .
$$ 
\end{lemm}
\begin{proof}
See Hoogenboom \cite[Theorem~6.10]{xhoo}
for the case $\tau \theta = \theta\tau$,
 and Matsuki \cite[ Theorem~1]{xmahgl} for the general case
where $\tau$ may not commute with $\theta$.
\end{proof}

Now, Lemma~\ref{lem:cptvisible} combined with 
Lemma~\ref{lem:7.5} implies Theorem~\ref{thm:Hvisible}
for a compact simple $G_U$.
Hence, we have shown 
Theorem~\ref{thm:Hvisible} in Case I.
\qed

\subsection{Proof of Theorem~\ref{thm:Hvisible} (compact case)}
\label{subsec:cpt2}
This subsection gives a proof of Theorem \ref{thm:Hvisible} in Case
II. 
Suppose that both $G_U / K_1$ and $G_U / K_2$ are
 compact Hermitian symmetric spaces.
We write $\theta_1$, $\theta_2$ for the corresponding involutive
 automorphisms of $G_U$.
Then,
 applying Lemma~\ref{lem:7.5} 
to $(\theta_1, \theta_2)$ in place of $(\theta,\tau)$,
 we find an involution $\sigma' \in \Aut(G_U)$
 satisfying the following three conditions:
\begin{align}
& \text{$\sigma' \theta_i = \theta_i \sigma'$ $(i =1, 2)$.}
\label{eqn:7.8.1}
\\
& \text{The induced action of $\sigma'$ on $G_U/K_i$ $(i=1,2)$ is anti-holomorphic.}
\label{eqn:7.8.2}
\\
& \text{$(\mathfrak{g}_U)^{\sigma', -\theta_1, -\theta_2}$
 contains a maximal abelian subspace of $(\mathfrak{g}_U)^{-\theta_1, -\theta_2}$.}
\label{eqn:7.8.3}
\end{align}
We remark that the condition \eqref{eqn:7.8.2} 
for $i=2$ is not included in Lemma~\ref{lem:7.5}, 
 but follows automatically by our choice of $\sigma'$.

Now, we define three involutive automorphisms $\tau$, $\theta$ and $\sigma$ on
 $G_U \times G_U$ by 
$\tau(g_1, g_2) := (g_2, g_1)$,
 $\theta:=(\theta_1, \theta_2)$
 and $\sigma := (\sigma', \sigma')$, respectively.
Then $(G_U \times G_U)^\tau = \diag(G_U)$.
By using the identification
$$
   (\mathfrak{g}_U \oplus \mathfrak{g}_U)^{-\tau}
   = \set{(X, -X)}{X \in \mathfrak{g}_U} 
   \rarrowsim \mathfrak{g}_U  ,
  \quad (X, -X) \mapsto X  ,
$$
we have
\begin{align*}
   (\mathfrak{g}_U \oplus \mathfrak{g}_U)^{-\tau, -\theta}
   &\simeq (\mathfrak{g}_U)^{-\theta_1, -\theta_2}  ,
\\
   (\mathfrak{g}_U \oplus \mathfrak{g}_U)^{\sigma, -\tau, -\theta}
   &\simeq (\mathfrak{g}_U)^{\sigma', -\theta_1, -\theta_2}  .
\end{align*}
Then it follows from \eqref{eqn:7.8.3} 
that
  $(\mathfrak{g}_U \oplus \mathfrak{g}_U)^{\sigma, -\tau, -\theta}$
 contains a maximal abelian subspace of
  $(\mathfrak{g}_U \oplus \mathfrak{g}_U)^{-\tau, -\theta}$.
Thus, we can apply Lemma~\ref{lem:cptvisible} to $\tau, \theta$ and
$\sigma \in \Aut (G_U \times G_U)$,
and therefore conclude that the diagonal action of $G_U$ on
$G_U / K_1 \times G_U / K_2$ is strongly visible.
Hence Theorem~\ref{thm:Hvisible} in Case II has been proved.
Together with Subsection \ref{subsec:cpt1},
we have proved Theorem \ref{thm:Hvisible} for compact case.

By Sections \ref{sec:3} and \ref{sec:4}, 
we have now completed the proof of Theorem \ref{thm:Hvisible}. 
\qed

\section{Visible actions of unipotent subgroups}
\label{sec:5}
This section gives a proof of the strong visibility of 
the maximal unipotent group $N$ action 
on the Hermitian symmetric space $G/K$ (Theorem~\ref{thm:Nvisible}).
The proof parallels to that for the $K$-action on $G/K$.

\begin{proof}[Proof of Theorem~\ref{thm:Nvisible}]
Without loss of generality,
we may and do assume that $G$ is a non-compact,
simply connected, simple Lie group of Hermitian type.
We take $\sigma$ to be the involution of $G$
as in Lemma~\ref{lem:4.1.1}.  

Let $\mathfrak{a}$ be a maximal abelian subspace in
$\mathfrak{g}^{\sigma,-\theta}$,
and set $A := \exp \mathfrak{a}$ and
$S := A \cdot o \subset G/K$.
We fix a positive system $\Sigma^+ (\mathfrak{a})$
of the restricted root system
$\Sigma(\mathfrak{a}) \equiv \Sigma(\mathfrak{g}, \mathfrak{a})$,
and define
$\mathfrak{n}_+ :=
 \sum_{\lambda \in \Sigma^+} \mathfrak{g} (\mathfrak{a}; \lambda)$.
Then we have an Iwasawa decomposition
$G = N_+ A K$ where
$N_+ = \exp (\mathfrak{n}_+)$
is a maximal unipotent subgroup of $G$.
Since $\mathfrak{a} \subset \mathfrak{g}^\sigma$,
$\sigma (\mathfrak{g} (\mathfrak{a} ; \lambda))
 = \mathfrak{g} (\mathfrak{a}; \lambda)$
for any $\lambda$.
In particular,
$\sigma$ stabilizes $\mathfrak{n}_+$.

Let $N_+ \cdot x$ be an $N_+$-orbit through $x \in G/K$.
We write $x = nak \cdot o$,
where $o = eK$.
Then $N_+ \cdot x = N_+ a \cdot o$,
and 
$\sigma(N_+ \cdot x) = \sigma(N_+) \sigma (a) \cdot o 
 = N_+ a \cdot o = N_+ \cdot x$.
Thus, $\sigma$ preserves each $N_+$-orbit on $G/K$.
Furthermore,
$\sigma$ acts anti-holomorphically on $G/K$ by
\eqref{eqn:2.4.3}
and $\sigma |_S = \mathrm{id}$
by $\mathfrak{a} \subset \mathfrak{g}^{\sigma}$.
Hence, the action of $N_+$ on $G/K$ is strongly visible.
Since any maximal unipotent subgroup $N$ is conjugate to $N_+$,
we have proved Theorem~\ref{thm:Nvisible}.
\end{proof}

\section{Applications to representation theory}
\label{sec:6}
In \cite{mfbdle},
we have given an abstract theorem on propagation of multiplicity-free
property of representations from fibers to spaces of sections for
equivariant holomorphic vector bundles.
Its main assumption is that actions on the base spaces are strongly
visible. 
Accordingly, if we find strongly visible actions on complex manifolds,
then we can expect a number of multiplicity-free theorems.

This section gives a brief explanation about how
 our geometric results 
(e.g.\ Theorems \ref{thm:Hvisible} and \ref{thm:Nvisible})
are applied to such multiplicity-free theorems by confining ourselves
 to the line bundle cases
(representations on fibers are automatically irreducible).
Detailed proof of these applications is given in a separate paper
\cite{mf-korea}.
Surprisingly, there was no literature,
to the best of our knowledge, before \cite{Saga1997},
on a systematic study of 
multiplicity-free theorems for the restriction with respect to general 
symmetric pairs,
although a number of explicit branching laws had been previously
known especially for finite dimensional representations
in the classical case
(see \cite{xhowemultone,xkratten,xokada} for example).
Our applications include both finite and infinite dimensional
representations,
and both discrete and continuous spectra. 

First, suppose $G$ is a connected Lie group of Hermitian type.
Retain the setting as in Subsection~\ref{subsec:3.3}.
Let $(\pi,\mathcal{H})$ be an irreducible unitary representation of
$G$, 
and $\mathcal{H}_K$ the underlying
$(\mathfrak{g}_{\mathbb{C}},K)$-module. 
Then, it is known that
the $K$-module
$
\mathcal{H}_K^{\mathfrak{p}_+} :=
\{v \in \mathcal{H}_K : d\pi(Y)v = 0
$
for any
$
Y \in \mathfrak{p}_+\}
$
is either zero or irreducible.
We say $(\pi,\mathcal{H})$ is an irreducible
\textit{unitary highest weight representation}
if $\mathcal{H}_K^{\mathfrak{p}_+} \ne \{ 0 \}$.
Furthermore, $\pi$ is of \textit{scalar type}
if $\dim \mathcal{H}_K^{\mathfrak{p}_+} = 1$.
Then, Theorem~\ref{thm:Hvisible} leads us to the following
multiplicity-free theorem:
\begin{corollary}[{\mdseries\cite{Saga1997}, 
\cite[Theorem A]{mf-korea}}\bfseries]
\label{cor:mfGH}
If $(\pi,\mathcal{H})$ 
is an irreducible unitary highest weight representation
of $G$ of scalar type,
then for any symmetric pair $(G,H)$,
the restriction $\pi |_H$ is multiplicity-free.
\end{corollary}

\begin{Rem}
Once we know the branching law is a priori multiplicity-free,
it would be interesting and reasonable to try to find its explicit formula.
It is noteworthy that ``new'' irreducible spherical unitary
representations of $H$ may occur as discrete summands in the setting
of Corollary \ref{cor:mfGH} (e.g.\ \cite{xorz}).
We remark that irreducible spherical unitary representations have not
been classified for general reductive Lie groups
(see \cite{xbarbasch}).
\end{Rem}

\begin{Rem}
\label{rem:mfGHbdle}
The multiplicity-free theorem for the vector bundle case
\cite{mfbdle} strengthens 
Corollary~\ref{cor:mfGH},
namely, one can relax the hypothesis of scalar type to
the following condition:
\addtocounter{subsection}{1}
\setcounter{equation}{0}
\begin{equation}
  \label{eqn:minKMF}
  \mathcal{H}^{\mathfrak{p}_+}_K \ 
  \text{is multiplicity-free as a $Z_{H\cap K}(\mathfrak{a})$-module.}
\end{equation}
Here, $Z_{H\cap K}(\mathfrak{a})$ 
is the centralizer of $\mathfrak{a}$ in $H \cap K$, and 
$\mathfrak{a}$ is a maximal abelian subspace of
$\mathfrak{g}^{-\theta, \sigma,\tau\theta}
 = \mathfrak{g}^{-\theta, \sigma,-\tau}$
if $H$ is defined by an involution $\tau$ 
and $\sigma$ is another involution given by 
Lemma \ref{lem:5.1}. 

The condition \eqref{eqn:minKMF} is obviously
satisfied if $\mathcal{H}^{\mathfrak{p}_+}_K$
is one dimensional (i.e. $\pi$ is of scalar type).
The smaller the dimension of the
totally real submanifold
$S = (\exp \mathfrak{a}) \cdot o$ is,
the larger the centralizer
 $Z_{H\cap K}(\mathfrak{a})$ becomes and the more likely the
condition \eqref{eqn:minKMF} is satisfied for the $K$-type 
$\mathcal{H}_K^{\mathfrak{p}_+}$.
This indicates how the slice $S$ plays a crucial role in the
multiplicity-free theorem.
\end{Rem}

\begin{corollary}
\label{cor:tensor1}
Suppose $\pi_1$ and $\pi_2$ are unitary highest weight modules of
scalar type.
Then, the tensor product representation $\pi_1 \otimes \pi_2$
is multiplicity-free.
\end{corollary}
The above tensor product is discretely decomposable.
On the other hand,
the following case corresponding to Theorem~\ref{thm:Gvisible}~(2)
contains continuous spectra in general:

\begin{corollary}
\label{cor:tensor2}
Retain the setting of Corollary~\ref{cor:tensor1}.
Then the tensor product $\pi_1 \otimes \pi_2^*$
is multiplicity-free.
\end{corollary}

The following is a descendent of Theorem~\ref{thm:Nvisible}.

\begin{corollary}[{\mdseries\cite[Theorem 33]{RIMS}}\bfseries]
\label{cor:6.6}
Let $G$ be a non-compact simple Lie group of Hermitian type,
and $N$ a maximal unipotent subgroup.
If $\pi$ is an irreducible unitary highest weight representation of
scalar type,
then the restriction $\pi |_N$ is multiplicity-free.
\end{corollary}

One of the simplest examples for Corollary \ref{cor:6.6}
 is the fact that the Hardy space
(an irreducible representation of
$G = SL(2,\mathbb{R})$
has simple spectra supported on the half line on the Fourier
transform side
(decomposition into irreducible representations of
$N \simeq \mathbb{R}$).

So far, we have discussed a non-compact $G$.
Next, let us consider the compact case.
Suppose now $G$ is a connected, compact simple Lie
group.
We take a Cartan subalgebra $\mathfrak{t}$ of $\mathfrak{g}$,
and fix a positive system
$\Delta^+ (\mathfrak{g}_{\mathbb{C}}, \mathfrak{t}_{\mathbb{C}})$.
For a dominant integral weight 
$\lambda \in \mathfrak{t}_{\mathbb{C}}^*$,
we write $\pi_\lambda$ for the irreducible finite dimensional
representation of $\mathfrak{g}_{\mathbb{C}}$
with highest weight $\lambda$.
We say that $\pi_\lambda$ is 
\textit{pan type} if $(\mathfrak{g},\mathfrak{g}_\lambda)$
forms a symmetric pair where
$\mathfrak{g}_\lambda$ is the Lie algebra of the
isotropy group 
$G_\lambda = \{g \in G: \textrm{Ad}^*(g)\lambda = \lambda\}$
(here, we regard 
$\mathfrak{t}^*_{\mathbb{C}} \subset \mathfrak{g}^*_{\mathbb{C}}$
via the Killing form), that is,
$$
\mathfrak{g}_\lambda :=
\{ Y \in \mathfrak{g} : \lambda ([Y,Z]) = 0
\text{ \ for any $Z \in \mathfrak{g} \}$}.
$$
See Richardson--R\"{o}hrle--Steinberg \cite{xrrs} or
\cite[Lemma 6.2.2]{RIMS} for the list of such
$\pi_\lambda$.

Then, the following two corollaries are obtained again from
Theorem~\ref{thm:Hvisible}. 
They may be regarded as finite dimensional versions of
Corollary~\ref{cor:mfGH}, 
and Corollaries~\ref{cor:tensor1} and \ref{cor:tensor2},
respectively.

\begin{corollary}
[{\mdseries\cite[Theorem~26]{RIMS}}\bfseries]
\label{cor:mfcptsymm}
Let $\pi_\lambda$ be a pan representation of a connected
compact Lie group $G$.
Then the restriction $\pi_\lambda |_H$ is multiplicity-free
with respect to any symmetric pair $(G,H)$.
\end{corollary}

\begin{corollary}[{\mdseries\cite[Theorem~25]{RIMS}, \cite{xlittel}}\bfseries]
\label{cor:mfcpttensor}
The tensor product representation $\pi_\lambda \otimes \pi_\mu$ of any
two pan representations $\pi_\lambda$ and $\pi_\mu$ is
multiplicity-free.
\end{corollary}

\begin{Rem}
\label{rem:tensornear}
As in Remark~\ref{rem:mfGHbdle}, we can strengthen 
Corollaries~\ref{cor:mfcptsymm} and \ref{cor:mfcpttensor}
by using the vector bundle case \cite{mfbdle}.
This strengthened version covers,
for instance, Stembridge's classification \cite{xstemb} of the pairs
$(\pi_1, \pi_2)$ 
of two irreducible finite dimensional representations of $GL_n$
such that $\pi_1 \otimes \pi_2$ is multiplicity-free.
It also covers some further multiplicity-free theorem such as
 the restriction of `nearly rectangular shape'
representations to symmetric subgroups (\cite{xkratten}).
\end{Rem}

Corollary~\ref{cor:mfcptsymm} could be proved alternatively
by a traditional approach,
that is, by showing the $H_{\mathbb C}$-sphericity
(the existence of an open orbit of a Borel subgroup):

\begin{prop}  \label{lem:Hspherical}
  Let $X$ be a compact Hermitian symmetric space,
and $G_{\mathbb C}$ the group of biholomorphic
transformations of $X$. Then $X$ is $H_{\mathbb C}$-spherical for
any $H_{\mathbb C}$ such that 
$(G_{\mathbb C}, H_{\mathbb C})$ is a complex symmetric pair.
\end{prop}

Conversely, 
Proposition~\ref{lem:Hspherical} 
is obtained as a corollary of Theorem~\ref{thm:Hvisible}
via Corollary~\ref{cor:mfcptsymm} and \cite{xvk}
(see \cite[Corollary~15]{RIMS}) from our view point.  

\begin{Rem}  \label{rem:nonpan}
Corollary~\ref{cor:mfcptsymm} holds for \textit{any}
symmetric pair $(G,H)$.
For individual pairs $(G,H)$, there are sometimes more families of
representations $\pi_\lambda$ such that
 the restrictions $\pi_{k\lambda}|_H$ are
 multiplicity-free for all $k \in \mathbb{N}$ even if they are
not pan representations, or equivalently, even if
$(\mathfrak{g},\mathfrak{g}_\lambda)$ are not symmetric pairs.

For example, 
for $(G,H) = (U(n), U(p) \times U(q))$ $(p + q = n)$,
this is the case if
\[
\mathfrak{g}_\lambda =
\begin{cases}
 \mathfrak{u}(1)^n &\text{if $\min(p,q) =1$,} \\
 \mathfrak{u}(n_1) + \mathfrak{u}(n_2) + \mathfrak{u}(n_3)
               &\text{if $\min(p,q) =2$,} \\
 \mathfrak{u}(n_1) + \mathfrak{u}(n_2) + \mathfrak{u}(n_3)\ (\min(n_1,n_2,n_3)=1)
               &\text{if $\min(p,q) \ge 3$,}
\end{cases}
\]
where $n = n_1+n_2+n_3$ (see \cite[Theorem~3.3]{Acta2004}).
See \cite{xkgencar} for an explicit construction of a totally real
submanifold $S$ that meets every $H$-orbit on
 the coadjoint orbit $\textrm{Ad}^*(G)\cdot \lambda$
(generalized flag variety).
We note that unlike the symmetric case treated in this article,
$S$ is not a flat submanifold in the
 non-symmetric $(\mathfrak{g},\mathfrak{g}_\lambda)$.
\end{Rem}

\section{Coisotropic, polar, and visible actions}
\label{sec:cpv}
We conclude this paper with some comments on the following three related
concepts on group actions:
\begin{enumerate}
    \renewcommand{\labelenumi}{\upshape \theenumi)}
\item 
(Strongly) visible actions on a complex manifold
(Definition \ref{def:S}).
\item 
Coisotropic actions on a symplectic manifold
(\cite{GS, HW}).
\item 
Polar actions on a Riemannian manifold
(\cite{PT, PT2}).
\end{enumerate}
Suppose a compact Lie group $H$ acts
on a symplectic manifold $M$ by symplectic automorphisms.
The action is called \textit{coisotropic} if one and hence all
principal $H$-orbits $H \cdot x$ are coisotropic with respect to the
symplectic form $\omega$, i.e.\
$T_x(H\cdot x)^{\perp\omega} \subset T_x(H\cdot x)$.
By \cite[Theorems 4 and 7]{RIMS},
Theorems \ref{thm:Hvisible} and \ref{thm:Nvisible} imply:
\begin{corollary}\label{cor:coisotropic}
{\upshape 1)}\enspace
In Theorem \ref{thm:Hvisible},
the $H$-action on $D$ is also coisotropic and visible.

{\upshape 2)}\enspace
In Theorem \ref{thm:Nvisible},
the $N$-action on $D$ is also coisotropic and visible.
\end{corollary}

On the other hand,
an isometric action of a compact Lie group $H$ on a Riemannian
manifold $M$ is called \textit{polar} if there exists a closed
connected submanifold $S$ of $M$ that meets every $H$-orbit
orthogonally.
By a classic theorem of R. Hermann \cite{xHermann},
the $H$-action on $G/K$ is polar if $G$ is compact and both $(G,H)$
and $(G,K)$ are symmetric pairs.

Since the polar action on an irreducible compact
homogeneous K\"{a}hler
manifold is coisotropic by \cite[Theorem 1.1]{PT} and visible by
\cite[Theorem 6]{RIMS}, Corollary \ref{cor:coisotropic} (1) for compact
$D$ follows also from Hermann's theorem \cite{xHermann}.
However, what we needed for strong visibility was not only the construction of 
a slice but also that of an anti-holomorphic diffeomorphism $\sigma$.
This was a core of the proof of Theorem
\ref{thm:Hvisible} in Sections \ref{sec:3} and \ref{sec:4}. 
We note that (strongly) visible actions are not always polar in
the non-symmetric case in general (see \cite{xkgencar}).

\end{document}